\documentclass[11pt]{article}

\usepackage[a4paper,margin=2.5cm]{geometry}

\usepackage{lipsum}
\usepackage{amsfonts}
\usepackage{graphicx}
\usepackage{epstopdf}
\usepackage{algorithm}
\usepackage{algorithmic}
\usepackage{amsmath,amssymb,amsfonts,amsthm}
\usepackage{amsopn}
\usepackage{bbm}
\usepackage{mathrsfs}
\usepackage{hyperref}
\usepackage[nameinlink,noabbrev]{cleveref}
\usepackage{xcolor}
\usepackage[normalem]{ulem}
\usepackage{comment}
\theoremstyle{plain}
\newtheorem{theorem}{Theorem}[section]
\newtheorem{lemma}[theorem]{Lemma}

\newtheorem{corollary}[theorem]{Corollary}
\theoremstyle{definition}
\newtheorem{definition}[theorem]{Definition}
\newtheorem{assumptions}[theorem]{Assumptions}
\theoremstyle{remark}
\newtheorem{remark}[theorem]{Remark}

\numberwithin{equation}{section}

\newenvironment{keywords}{\par\smallskip\noindent\textbf{Keywords. }}{\par\smallskip}
\newenvironment{MSCcodes}{\par\smallskip\noindent\textbf{MSC Classification. }}{\par\smallskip}
\newcommand{\email}[1]{\texttt{#1}}

\title{Infinite-Horizon Optimal Control of Jump-Diffusion Models for Pollution-Dependent Disasters}
\author{Daria Sakhanda\thanks{ETH Z\"urich, Department of Mathematics, Switzerland (\email{daria.sakhanda@math.ethz.ch}).}
	\and Joshué Helí Ricalde-Guerrero\thanks{ETH Z\"urich, Department of Mathematics, Switzerland (\email{joshue.ricalde@math.ethz.ch}).}}

\def\cA{{\cal A}}

\def\cF{{\cal F}}

\def\cL{{\cal L}}
\def\cS{{\cal S}}
\def\cT{{\cal T}}

\def\pt{\partial}
\def\n{\noindent}

\makeatletter

\providecolor{lyxadded}{rgb}{0,0,1}
\providecolor{lyxdeleted}{rgb}{1,0,0}

\newcommand{\veps}{\varepsilon}

\newcommand{\bbE}{\ensuremath{\mathbb E}}
\newcommand{\bbF}{\ensuremath{\mathbb F}}

\newcommand{\bbP}{\ensuremath{\mathbb P}}

\newcommand{\bbR}{\ensuremath{\mathbb R}}

\begin{document}
	\date{}
	\maketitle
	
	\begin{abstract}
		This paper is devoted to developing a unified framework for stochastic growth models with environmental risk, in which rare but catastrophic shocks interact with capital accumulation and pollution. 
		The analysis is based upon a general Poisson point process formulation, leading to non-local Hamilton-Jacobi-Bellman (HJB) equations that admit closed-form candidate solutions and yield a composite state variable capturing exposure to rare shocks.
		We consider cases where disaster risk is endogenized through a pollution-dependent intensity and, in the more general cases, it also accommodates for state-dependent events of varying magnitude. Our formulation captures how environmental degradation amplifies macroeconomic vulnerability and strengthens incentives for abatement. 
		From a technical perspective, it provides tractable jump-diffusion control problems whose HJB equation decomposes naturally into capital and pollution components under power-type value function. 
	\end{abstract}
	


\section{Introduction}

\n 

The study of economic growth under environmental risk raises fundamental  mathematical questions in stochastic control for systems combining continuous  dynamics and discontinuous jumps. Classical growth-environment models typically describe pollution and damages as deterministic or smoothly evolving, with dynamics  governed by stochastic differential equations without discontinuities. 
Empirical evidence, however, shows that risk also takes the form of rare but catastrophic events -- such as abrupt climate disasters, ecosystem collapses, or large-scale technological failures -- which arrive unpredictably and are naturally modeled as jumps with large economic impacts.

From the mathematical side, the challenge is to develop a rigorous and tractable framework that  incorporates both continuous fluctuations (e.g.\ Brownian pollution dynamics) and state-dependent jump risks (e.g.\ disaster intensities  increasing with pollution). This necessitates the use of results from the theory of optimal control such as the Dynamic Programming Principle (DPP), Hamilton-Jacobi-Bellman (HJB) equations and verification theorems that, unlike with their deterministic/diffusive counterparts, require the careful manipulation of non-local terms.

Motivated by economic questions of economic growth, environmental change, and climate risk, this paper aims to develop a unified framework for stochastic growth models with environmental risk, in which rare but catastrophic shocks interact with capital accumulation and pollution. More precisely, this paper extends the model framework of \cite{brausmann_escaping_2024} and contributes to the literature on \textit{growth under environmental risk} by organizing and progressively generalizing the stochastic framework for growth-environment problems:
\begin{enumerate}
	\item Starting from the benchmark case of a \emph{homogeneous Poisson process} where a social planner is tasked with allocating resources in an economy subject to homogeneous climate-related disasters -- occurring at a constant rate -- that destroy a fraction of the capital stock.
    
	\item Stepping up from the benchmark case, we analyze the case of an economy where the rate of occurrence of disasters is directly dependent on the pollution stock. That is to say, the destruction of capital is driven by a \emph{nonhomogeneous Poisson process}, representing an intrinsic feedback from environmental degradation to disaster risk.  
   
	\item Next, we introduce \emph{Brownian noise in the pollution dynamics} to better represent the possible continuous fluctuations within the system, as opposed to only discrete random shocks.
  
	\item Finally, we generalize the entire framework using \textit{Poisson random measures with marks}, in order to allow for disasters of varying magnitudes and more general jump structures (e.g. small vs. large forest fires in an already deforested community).
   
\end{enumerate}

To the best of our knowledge, models 2--4 are new to the literature at the time of writing (the benchmark case 1 has been addressed in \cite{brausmann_escaping_2024}). Furthermore, these models balance generality and tractability, providing a rigorous analytical foundation as well as practical tools for analyzing and solving control problems that govern the long-run interaction between the economy and an environment with rare but potentially catastrophic events. Thus, this paper contributes both to the applied theory of sustainable growth and the broader mathematical toolbox for stochastic control with jumps. 

On the modelling side, we organize a broad class of pollution-driven disaster specifications -- ranging from constant-intensity shocks to marked Poisson random measures -- within a unified stochastic control framework. This clarifies the relationship between existing models and how richer forms of environmental risk can be incorporated without sacrificing tractability. 

On the analytical side, we explicitly derive the associated HJB equation for each model and establish verification results that rigorously characterize the planner’s value function under general assumptions. We do this by exploiting the relation between \textit{state} and \textit{costate} variables present in the Hamiltonian of the DPP, in order to carefully construct a candidate for the value function of the social planner's problem. 

The rest of the article is organized as follows. In Section \ref{section:Models_Framework}, we specify all of the elements in our setting and formally introduce the models under study. Section \ref{Section:Hamilton-Jacobi-Bellman-PIDEs} is devoted to the technical analysis of our framework based on HJB Partial Integro-Differential Equations (PIDEs): in the first half -- Section 3.1 -- we review the essential tools and techniques for our study, while on the second half -- Section 3.2 -- we examine the core components that link our models together, as well as each of the individual HJB equations. For the sake of readability, we shall leave all the corresponding technical computations of this section to Appendix \ref{Sec:Technical-appendix}. Lastly, Section \ref{Sec:Verification-theorems-viscosity-solutions} is dedicated to our main results: a collection of model-specific verification theorems addressing state-dependent disasters and diffusion-driven uncertainty. Particularly, for model 4 we provide a forward-backward stochastic differential equation representation of the value function, linking the control problem to FBSDE methods.

\begin{remark}
Throughout the paper we will mostly follow the notation used in \cite{brausmann_escaping_2024}. Moreover, unless otherwise explicitly stated, we differentiate between the random values of a stochastic process and one of its realizations (i.e. an arbitrary deterministic value) by the appearance of the time subindex, e.g. between $X_t$ and $X$, respectively.
\end{remark}

\section{Modelling Framework}
\label{section:Models_Framework}

\n 

We consider a stylized representation of the global economy, which produces a single composite good under constant returns to scale. Production relies on the aggregate capital stock at time $t$, denoted by $K_t$ which includes physical capital, human capital, and intangible assets. The production process $(K_t)_{t \geq 0}$ generates pollution: at each instant $t$, greenhouse gas (GHG) emissions $E_t$ are released into the atmosphere. These emissions accumulate in the atmospheric pollution stock, $P_t$, which increases with the flow $E_t$ and decreases at a natural absorption rate $\alpha \in [0,1)$, assumed to be small or negligible. While $P_t$ is referred to as the pollution stock, it can more generally be interpreted as the \textit{inverse of environmental quality}. Similarly, the emissions variable $E_t$, can be viewed more broadly as any environmentally damaging by-product of economic activity.

The model incorporates the possibility of natural disasters (we refer to it as an ``event'') which occur randomly over time. When such an event occurs, it instantaneously destroys a fraction of the capital stock, with the \textit{surviving share of capital} given by 
\begin{align*}
    \omega: \mathcal{S} \to(0,1),
    \qquad
    \cS := \bbR_{>0} \times \bbR_{> 0},
\end{align*}
which is determined endogenously and depends on the current levels of pollution and capital:
\begin{equation}
	\omega(K, P) = e^{- \delta P^{\xi} K^{\eta}},
    \qquad 
    \forall (K,P) \in \mathcal{S},
\end{equation}
where $\xi \geq 0$ and $\eta \geq 0$ are parameters that capture the sensitivity of damage intensity to pollution and capital, respectively. 

The \textit{output of production}, 
\begin{align*}
    Y: \bbR_{> 0} \to \bbR_{\geq 0},
\end{align*}
can be allocated to consumption, investment in capital, or environmental protection. Following \cite{brausmann_escaping_2024}, we set a constant rate of production $A > 0$ so that 
\begin{align*}
    Y(K) := AK,
    \qquad \forall K > 0.
\end{align*}

A fraction $\theta\in[0, 1]$ of output is allocated to \textit{abatement}, yielding \textit{abatement investment} $I: [0, 1] \times \bbR_{>0} \to \bbR_{\geq 0}$, defined as
\begin{equation}
    I(\theta, K) = \theta Y (K),
    \qquad 
    \forall (\theta , K) \in [0,1] \times \mathbb{R}_{> 0}.
\end{equation}
The \textit{remaining share}, $(1 - \theta) Y$, is split between consumption and capital accumulation. 

\textit{Abatement activities} reduce emissions through a monotone increasing function $Z: \bbR_{\geq 0} \to \bbR_{\geq 0}$, defined as 
\begin{equation} \label{eq:Z_I}
	Z(I) = z I,
\end{equation}
where $z > 0$ stands for the efficiency of abatement. In economical terms, abatement activities exhibit a constant return $z$, with the \textit{total abatement} $Z$ proportional to the resources allocated to it.  

On the other hand, we assume the \textit{total emissions} 
$$E: \bbR_{\geq 0 }\times \bbR_{>0} \to \bbR_{\geq 0}$$ 
are \textit{homogeneous in time}; that is, their dependence on time occurs only through the current state of the system. More precisely, it is given by the net balance between \textit{gross emissions} and the \textit{total abatement}:
\begin{equation*} 
    E(Y,Z) = \phi Y - Z,
\end{equation*}
where  $\phi > 0$ is a proportionality constant describing the intensity of emissions due to production. Equivalently, in terms of \textit{abatement investment} and \textit{capital}: 
\begin{equation} 
    \label{eq:Model_Emission}
    E(\theta,K) = \phi Y(K) - \sigma I(\theta,K),
    \qquad
    \forall (\theta , K) \in [0,1] \times \mathbb{R}_{> 0}.
\end{equation}

The \textit{instantaneous preferences} of the planner are given by a twice continuously differentiable power-type utility function $U \in C^2(\cS)$, defined as
\begin{align}
    \label{eq:Utility_Fn}
    U: \bbR_{>0} \times \bbR_{> 0} \to \bbR,
    \qquad
    U(C, P) 
    := 
    \frac{C^{1 - \varepsilon}}{1 - \varepsilon} 
        - 
        \chi \frac{P^{1 + \beta}}{1 + \beta},
\end{align}
where
\begin{itemize}
    \item 
        $C > 0$ represents the \textit{allocated consumption};
    \item 
        $\varepsilon \in \mathbb{R}_{> 0} \setminus \{ 1 \}$ measures \textit{relative risk aversion};

    \item 
        $\beta > 0$ governs the \textit{curvature of disutility from pollution};

    \item 
        $\chi > 0$ reflects the \textit{weight placed on pollution} in the utility function.
\end{itemize}

To ensure the well-posedness of the optimization problem (that is, $U$ satisfies the regularity conditions required in dynamic optimization and Hamilton-Jacobi-Bellman frameworks, see \cite{dixit_investment_2012}), we verify standard properties of the utility function in \eqref{eq:Utility_Fn}:
\begin{enumerate}
    \item 
        \textit{Monotonicity and concavity.} For fixed $P$, the marginal utility of consumption is
        \begin{equation}
            \frac{\partial U}{\partial C}(C,P) = C^{-\varepsilon} > 0, 
            \qquad 
            \frac{\partial^2 U}{\partial C^2}(C,P) = -\varepsilon C^{-\varepsilon-1} < 0,
        \end{equation}
        for all $C>0$, implying that $U$ is strictly increasing and strictly concave in $C$. 
        
        Similarly, for a fixed $C$, we have
        \begin{equation}
            \frac{\partial U}{\partial P}(C,P) = -\chi P^{\beta} < 0,
            \qquad
            \frac{\partial^2 U}{\partial P^2}(C,P) = -\chi \beta P^{\beta-1} < 0,
        \end{equation}
        for all $P>0$, implying that $U$ is strictly decreasing and concave in the pollution stock $P$.
      
    \item 
        \textit{Inada conditions in consumption.} For all $P>0$, the marginal utility of consumption satisfies
        \begin{equation}
            \lim_{C \to 0^+} \frac{\partial U}{\partial C}(C,P) = +\infty,
            \qquad
            \lim_{C \to \infty} \frac{\partial U}{\partial C}(C,P) = 0.
        \end{equation}
        Thus, the Inada conditions hold with respect to $C$, ensuring the existence of interior optimal consumption paths under standard regularity assumptions.

    \item 
        \textit{Boundedness properties.} The map 
        \begin{align*}
            C \longmapsto \frac{C^{1 - \varepsilon}}{1 - \varepsilon}
        \end{align*}
        is strictly increasing and strictly concave on $\bbR_{> 0}$. Its behavior as $C \to \infty$ depends on the parameter $\varepsilon$ as follows:
        \begin{equation}
            \lim_{C\to\infty} \frac{C^{1-\varepsilon}}{1-\varepsilon}
            =
            \begin{cases}
                +\infty,    & \varepsilon \le 1,
                \\[4pt]
                0,          & \varepsilon > 1.
            \end{cases}
        \end{equation}
        Consequently, $U(C,P)$ is bounded above for each fixed $P>0$ if and only if $\varepsilon>1$. Moreover, for each fixed $C > 0$, the disutility from pollution satisfies 
        \begin{equation}
            \lim\limits_{P \to 0^{+}} -\chi \frac{P^{1 + \beta}}{1 + \beta} = 0 \quad \text{ and } \quad \lim\limits_{P \to \infty} -\chi \frac{P^{1 + \beta}}{1 + \beta} = -\infty.
        \end{equation}
As a result, the utility function $U(C, P)$ is unbounded below as $P \to \infty$. 
\end{enumerate}

Given $(K_0,P_0)\in \cS$, the social planner aims to maximize the expected discounted utility over an infinite horizon by choosing optimal paths for consumption $(C_t)_{t \geq 0}$ and abatement share $(\theta_t)_{t \geq 0}$,
\begin{equation}\label{eq:Optimal_v}
	v(K_0,P_0)=\sup_{(C_t,\theta_t)}\ \mathbb E \left[\int_0^\infty e^{-\rho t}U(C_t,P_t)\,dt\right],
\end{equation}
where the supremum is taken over all \textit{admissible processes} $(C_t,\theta_t)_{t \geq 0} \in \mathcal{A}$ (to be specified later) and $\rho > 0$ defines the constant rate of time preference. At the same time, the state process $(K_t, P_t)$ is subject to \textit{model-specific jump terms} and common drift components
\begin{align}
    \label{eq:b_K}
	b^{\mathrm{cap}} (K_t,P_t,C_t,\theta_t) &:=(1-\theta_t)Y(K_t)-C_t,
    \\
    \label{eq:b_P}
	b^{\mathrm{pol}}(K_t,P_t,C_t,\theta_t) &:=\phi Y(K_t)-Z(\theta_t Y(K_t))-\alpha P_t.
\end{align}
The planner's problem is then to \textit{optimally allocate resources between consumption and abatement in order to balance economic growth, environmental quality, and resilience against environmentally driven disasters}.

We analyze four increasingly rich variants of the model \eqref{eq:Optimal_v}, which share the state $(K_t,P_t)$, the controls $(C_t, \theta_t)$, the drift parts \eqref{eq:b_K}--\eqref{eq:b_P} and the preferences \eqref{eq:Utility_Fn}:
\begin{enumerate}
    \item 
        \textit{Constant arrival rate (Homogeneous Poisson Process, HPP).} Events arrive with constant intensity $\lambda>0$; jumps are unit-sized and destroy a state-dependent fraction of capital via $\omega(K_t,P_t)$.
	
    \item 
        \textit{Pollution-driven intensity (Nonhomogeneous Poisson Process, NHPP).} The intensity becomes state-dependent, $\lambda(P_{t-})$ (e.g., affine $\lambda_0 + \lambda_1 P_{t-}$), introducing feedback from environmental quality to disaster risk.
	
    \item 
        \textit{Jump-diffusion with stochastic pollution.} We keep $\lambda(P_{t-})$ and add Brownian fluctuations in $P_t$, yielding a diffusion term into the HJB.
        
	\item 
        \textit{Jump-diffusion with stochastic pollution and marked jumps (Poisson Random Measure, PRM)}. We generalize the previous model by incorporating a marked point process $q(dt,d\zeta)$ determined by the compensator $\lambda(P_{t-},\zeta)dt\otimes\nu(d\zeta)$. This change allows random disaster magnitudes to be considered in the model, therefore linking the problem to that of nonlocal operators and Partial Integro-Differential Equations (PIDEs).
\end{enumerate}

Since every model modifies the jump mechanism and/or the law of the pollution stock $(P_t)_{t \geq 0}$, it is necessary to work on a stochastic basis $(\Omega, \cF, \bbF, \bbP)$ fitted for each individual case. However, due to the similarity between the components, a unifying framework remains in place. We will address this in detail later, in Section \ref{Section:Hamilton-Jacobi-Bellman-PIDEs}.

\subsection{Baseline: Homogeneous Poisson Process (HPP)}
\label{section:Model under the constant arrival rate: standard Poisson Process}

\n 

On a filtered probability space $(\Omega,\mathcal F,\mathbb P)$, we define a Poisson process $q = (q_t)_{t \geq 0}$ with intensity $\lambda > 0$, that is, an $\bbF$-adapted c\`adl\`ag process with values in $\mathbb{N}$ such that 
\begin{enumerate}
\renewcommand{\theenumi}{\roman{enumi}}

    \item 
        $q_0 = 0$ almost surely, 

    \item 
        $q$ is continuous in probability, 

    \item 
        the increments are stationary and independent, i.e., for all $0 \leq s < t$, the random variable $q_t - q_s$ is independent of $\cF_s$ and $q_t-q_s\sim\mathrm{Poisson}\big(\lambda(t-s)\big)$.
\end{enumerate}

The parameter $\lambda$ can be characterized as
\begin{equation} \label{eq:lambda_det}
	\lambda := 
	\lim_{h \to 0} \frac{1}{h} \bbP[q_{t} - q_{t-h} = 1].
\end{equation} 
Over a small interval $[t-h,t]$, the probability of two or more jumps satisfies $\bbP \big[q_{t}-q_{t-h}\ge2\big]=o(h^2)$ as $h \to 0$, therefore we consider only scenarios with $0$ or $1$ jump in that interval. Then for $t > 0$, we can define the jump size of the Poisson process $q$,
\begin{equation} \label{Eq:Delta_q}
	\Delta q_h:= q_h - q_{h-} \in \{0, 1\}, \quad h \geq 0,	
\end{equation}
where $q_{t-} := \lim\limits_{s \nearrow t} q_s$, for all $t > 0$.	

In the model, each unit jump of $q$ represents the occurrence of a natural-disaster \emph{event}. Thus, events arrive according to a Poisson process with mean arrival rate $\lambda>0$. Given the initial $(K_0,P_0) \in \cS$, capital and pollution dynamics are given by
\begin{align}
	dK_t &= b^{\mathrm{cap}} (K_t,P_t,C_t,\theta_t)\,dt - \big(1-\omega(K_{t-},P_t)\big)\,K_{t-}\,dq_t,\label{Opt_K}
    \\
	dP_t &= b^{\mathrm{pol}} (K_t,P_t,C_t,\theta_t)\,dt. \label{Opt_P}
\end{align}

We next relax the constant-hazard assumption by letting the arrival rate respond to the contemporaneous pollution stock $P_t$.

\subsection{Pollution-Driven Disaster Intensity via Nonhomogeneous Poisson Process (NHPP)}
\label{section:Pollution-Driven-Disaster-Intensity-via-NHPP}

\n 

The baseline model in Section~\ref{section:Model under the constant arrival rate: standard Poisson Process} assumes that the expected frequency of destructive events is independent of the state of the environment. However, this assumption neglects the substantial empirical evidence linking environmental degradation to increased disaster likelihood. 

Rising concentrations of greenhouse gases, for instance, have been associated with heightened frequency and severity of extreme climatic events, such as storms, floods, and droughts.
To capture the empirically motivated idea that environmental degradation amplifies the frequency of natural disasters, we generalize the arrival process to a nonhomogeneous Poisson process whose intensity depends on the current stock of pollution $P_t$. Specifically, we replace the constant intensity $\lambda$ with a state-dependent stochastic intensity $\lambda(P_{t-})$, where
\begin{equation}
	\lambda: \bbR_{\geq 0} \to \bbR_{\geq 0}, \quad \lambda \in C^1, \quad \lambda'(P) > 0.
\end{equation}
Then, rare destructive events arrive following a simple point process $(q_t)_{t\ge0}$ such that, conditional on $\cF_{t_0}$ with $t_0 < t$, the increments satisfy
\begin{equation}
    q_t - q_s \sim \mathrm{Poisson} \left(\int_{t_0}^t \lambda(P_{s-})\,ds\right),
    \qquad 
    0\le t_0<t,
\end{equation}
so that $\Delta q_t\in\{0,1\}$ a.s.

We specify the pollution-dependent hazard rate as
\begin{equation} 
    \label{eq:NonHP_lambda}
    \lambda(P) 
    = 
    \lambda_0 + \lambda_1 P, \qquad \lambda_0 \geq 0, \, \lambda_1 \geq 0,
\end{equation}
where $\lambda_0$ represents the baseline hazard rate unrelated to environmental conditions, and $\lambda_1$ measures the marginal increase in disaster risk per unit of pollution stock. 

This specification ensures that hazard rates rise linearly with environmental damage, allowing the model to reflect the dual role of abatement: reducing both the direct disutility from pollution and the frequency of capital-destroying events. From a technical perspective, the affine form of $\lambda(\cdot)$ preserves much of the tractability of the homogeneous case while introducing a meaningful state dependence in the jump intensity, thereby enriching the policy implications without introducing the analytical complexity of fully general L\'evy jump structures.

In particular, the map $\lambda: \bbR_{\geq 0} \to \bbR_{\geq 0}$ is continuous, locally Lipschitz, and of at most linear growth, ensuring
integrability of the compensator and well-posedness of the jump term. Thus, at each jump time of $q$, a fraction $\omega(K_{t^-}, P_{t^-})\in(0,1)$ of the capital stock survives, i.e.,
\begin{equation}
	K_t = \omega(K_{t-}, P_{t-}) \, K_{t-}\quad \text{if} \quad \Delta q_t=1,
\end{equation}
with pollution unchanged contemporaneously, $P_t=P_{t-}$.

The dynamics of capital and pollution then become
\begin{align}
	dK_t &= b^{\mathrm{cap}} ( K_t,P_t,C_t,\theta_t )\,dt - \big(1-\omega(K_{t-},P_t)\big)\,K_{t-}\, d q_t,\label{eq:Opt_K_NHPP}
    \\
	dP_t &= b^{\mathrm{pol}} ( K_t,P_t,C_t,\theta_t )\,dt. \label{eq:Opt_P_NHPP}
\end{align}

Heuristically, the model proposed here describes a \textit{weak form of self-exciting interaction} in which capital and pollution are coupled through a bidirectional feedback, linking a degrading environment with economic vulnerability. Indeed, although lower abatement investment clearly leads to environmental degradation, the model further implies that an increase in the incidence of environmental disasters correspondingly amplifies capital destruction.

In contrast to the baseline model where the controller was fully subjected to the exogenous timing of disasters, the present model allows control actions to delay the occurrence of adverse events that reduce capital. This coupling between the action of the planner and the occurrence of disasters is only strengthened in the models to come by further specifying how environmental harm amplifies economic losses.

\subsection{Jump-Diffusion Pollution with Intensity Feedback}
\label{section:BM+NHPP-driven-model}

\n 

We now extend our framework by allowing the pollution stock itself to be subject to stochastic fluctuations driven by Brownian motion, while retaining the specification of rare destructive events through a nonhomogeneous Poisson process whose intensity depends on the current level of pollution. This enriches the dynamics of the economy by capturing both state-dependent disaster risk and continuous environmental uncertainty.

The dynamics of capital $K_t$ remain as in \eqref{eq:Opt_K_NHPP}, while the dynamics of pollution $P_t$ with initial condition $P_0 > 0$ are given by
\begin{equation}\label{eq:P_W_NonH_Pois}
	dP_t = b^{\mathrm{pol}} ( K_t,P_t,C_t,\theta_t )\, dt + \sigma P_t \,dW_t,
\end{equation}
where $W = (W_t)_{t \geq 0}$ is a Brownian motion and $\sigma > 0$ is the diffusion parameter. The multiplicative term $\sigma P_t\,dW_t$ captures proportional (log-normal-type) fluctuations in pollution (e.g., meteorological dispersion or natural absorption shocks) and ensures strict positivity of $P_t$ for $P_0>0$, so that the state space $(K_t,P_t)\in\mathcal S$ is preserved without the need for boundary conditions at $P=0$.

Because the disaster intensity $\lambda(P_t)$ is increasing in $P_t$, stochastic pollution generates time-varying and endogenous disaster exposure. Abatement therefore plays a dual role: it reduces both the level of pollution damages and the level and variability of disaster risk. In particular, the conditional variance satisfies
\begin{equation}
\text{Var}(dP_t \mid \mathcal{F}_t) = \sigma^2 P_t^2\,dt,
\end{equation}
so absolute volatility rises with pollution, while relative volatility remains constant, consistent with proportional environmental uncertainty at the aggregate level.

From a technical perspective, \eqref{eq:P_W_NonH_Pois} leads to an integro-diffusion HJB equation with second-order term (see \eqref{Eq:2nd-order-HJB-term} below),
placing the problem within the standard class of stochastic control models with positive diffusions (e.g., \cite{dixit_investment_2012}). The deterministic benchmark is recovered as $\sigma \to 0$. The geometric diffusion in \eqref{eq:P_W_NonH_Pois} therefore provides a tractable and economically natural specification that preserves positivity and captures level-dependent environmental risk.

\subsection{Generalized Framework with Poisson Random Measures (PRMs)}
\label{Eq:PRM-Model-Introduction}

\n 

The Poisson random measure (PRM) formulation provides a unified and rigorous framework for modeling jump processes by making disaster arrivals and their compensators explicit and decomposing capital dynamics into drift and martingale components. This generalization not only strengthens the mathematical foundations of the model, but also facilitates extensions, such as allowing for disaster magnitudes. In this way, the PRM formulation serves as a bridge between the models introduced in Sections~\ref{section:Model under the constant arrival rate: standard Poisson Process}--\ref{section:BM+NHPP-driven-model} and richer specifications with marked jumps. 

Within this PRM framework, the social planner chooses consumption and abatement policies to balance output growth, environmental quality, and resilience to pollution-driven disasters, whose arrivals are governed by a state-dependent intensity (i.e. the frequency at which disasters occur) and whose magnitudes may be modelled as random marks. 
For example, within our setting it is possible to consider a more realistic economy in which

\begin{itemize}
    \item 
        the higher the pollution stock $P_t$, the more frequent disasters become;
    \item 
        the higher the pollution stock $P_t$, the more destructive disasters become.
\end{itemize}
Both channels amplify capital losses and environmental risk.

Mathematically, we consider a Poisson random measure 
\begin{align*}
    q : (\Omega,\mathcal{F},\mathbb{P}) \longrightarrow \mathcal{M}^*_c \big( [0,\infty) \times [0,\infty) \big),
\end{align*}
where $(\Omega, \mathcal{F}, \mathbb{F}, \mathbb{P})$ is a prescribed probability space and $\mathcal{M}^*_c( \mathcal{X} )$ denotes the space of simple, counting Borel measures over $\mathcal{X}$, endowed with the topology of weak convergence, see \cite{daley_introduction_2003}--\cite{daley_introduction_2008}. This means that $q$ is of the form
\begin{align}
    \label{Eq:q-proper}
    q(d t, d \zeta)
    =
    \sum_{n \geq 1} \delta_{ (\tau_n, \Delta_n) }(d t, d \zeta),
\end{align}
where $\tau_n$ denotes the occurrence time of the $n$-th disaster, and $\Delta_n$ its magnitude. In fact, we have that $\Delta_n$ follows the conditional law
\begin{align}
    \label{Eq:Conditional-mark-distribution}
    \mathbb{P} \big[ \Delta_n \in B \big\vert \sigma \{ (\tau_1, \Delta_1), \ldots, (\tau_{n-1}, \Delta_{n-1} ) \} \vee \sigma\{ \tau_n \} \big]
    \propto
    \int_{ B } \lambda( P_{\tau_n}, \zeta ) \nu( d \zeta ),
\end{align}
for all $B \in \mathcal{B}( \mathbb{R}_{\geq 0} )$, up to a normalizing constant. 

Regarding the compensator of $q$, it is assumed to have the form
\begin{align}
    \label{Eq:lambda}
    \Lambda( dt, d\zeta ) := \lambda(P_{t},\zeta) d t \otimes \nu(d \zeta),
\end{align}
where $\nu$ is a finite measure over $\mathbb{R}_{ \geq 0 }$ with finite second moment, and 
\begin{align*}
    \lambda: \mathbb{R}_{\geq 0} \times \mathbb{R}_{\geq 0} \to \mathbb{R}_{\geq 0}
\end{align*}
is a continuous differentiable function with
\begin{align*}
    &\frac{\partial}{\partial P} \lambda(P,\zeta) > 0, \qquad (P,\zeta) \in \mathbb{R}_{>0} \times \mathbb{R}_{\geq 0}
\end{align*}
and
\begin{align*}
    &\int_{ [0,\infty) } (1 + \zeta^2) \lambda(P,\zeta) \nu(d \zeta) < \infty,
    \qquad \forall P > 0.
\end{align*}

\begin{remark}
\label{Remark:Random-Intensity}
In order to allow models in which the disaster arrival rate is endogenous (i.e., the intensity $\lambda$ of the jump measure $q$ may depend on $P$), $q$ has to be a random measure able to sustain a (possibly) random intensity measure $\Lambda$. To overcome this technicality, it is necessary to ensure our stochastic basis $(\Omega, \mathcal{F},\mathbb{F}, \mathbb{P})$ is sufficiently rich. We will address this issue more precisely at the beginning of the next section.    
\end{remark}

With these assumptions in place, the controlled dynamics now become
\begin{align}
    \label{Eq:K-BM}  
	dK_t 
        &= 
        b^{\mathrm{cap}}(K_t,P_t,C_t,\theta_t)\,dt 
        - 
        \int_{ [0,\infty) } \big( 1 - \omega( K_{t-}, P_t, \zeta ) \big) K_{t-}\,q( dt \times d\zeta),
    \\
    \label{Eq:P-BM}
    dP_t 
        &= 
        b^{\mathrm{pol}}(K_t,P_t,C_t,\theta_t)\,dt + \sigma P_t\,d W_t,
    \\
    \label{Eq:K_0-P_0-BM}
    K_0 & > 0,\quad P_0  > 0,
\end{align}
where $b^{\mathrm{cap}}$ and $b^{\mathrm{pol}}$ are as in \eqref{eq:b_K}--\eqref{eq:b_P}, $( C_t )_{ t \geq 0 }$ and $( \theta_t )_{ \geq 0 }$ are $\mathbb{F}$-adapted processes such that the system \eqref{Eq:K-BM}--\eqref{Eq:K_0-P_0-BM} is well-defined, and
\begin{align*}
    (K,P,\zeta) \longmapsto \omega(K,P,\zeta) \in (0,1)
\end{align*}
represents the surviving proportion of capital after a disaster, which is dependent on the current levels of capital $K$ and pollution $P$, as well as the magnitude of the disaster $\zeta$. 

The function $\omega$ is assumed to be continuous in all its coordinates, and to satisfy the integrability condition
\begin{align*}
    &\int_{ [0,\infty) } \omega(K,P,\zeta)^2 \lambda(P,\zeta) \nu(d \zeta) < \infty \quad \forall (K,P) \in \cS, \quad \forall t \geq 0.
\end{align*}

Moreover, we assume $\omega$ is decreasing on $K$ (resp. $P$ and $\zeta$), thereby reflecting a higher vulnerability when capital intensity (resp. level of pollution and magnitude of the disaster) increases.

As stated previously, this model is motivated by the notion that the level of pollution affects not only the frequency but also the magnitude of the risk involved. 
For example, consider the following: let $q$ be a Poisson random measure determined by its compensator
\begin{align}
    \label{Eq:Compensator-Example}
    \Lambda( dt, d\zeta ) 
    :=
    \frac{\lambda( P_{t-} )}{\Gamma(P_{t-})} \zeta^{P_{t-} - 1} e^{-\zeta}\, dt \otimes d\zeta,
\end{align}
where $\lambda : \mathbb{R}_{\geq 0} \to \mathbb{R}_{\geq 0}$ is a measurable function; that is, the kernel $\Lambda$ corresponds to $\lambda(P_{t-})$ times the density of a Gamma distribution of shape parameter $P_t$ and scale parameter $1$. 

As a result, from \eqref{Eq:q-proper} and \eqref{Eq:Conditional-mark-distribution} we have that the disasters $q$ arrive as a Poisson point process of intensity 
\begin{align*}
    \mathbb{E}\Big[ q \big( (0,{t}] \times \mathbb{R}_{ \geq 0 } \big) \Big] 
    = 
    \int_{(0,t]} \lambda(P_{t-}) \Big( \int_{ \mathbb{R}_{ \geq 0 } } \frac{1}{\Gamma(P_{t-})} \zeta^{P_{t-} - 1} e^{-\zeta}\, d\zeta\Big)\, dt
    =
    \int_{(0,t]} \lambda(P_{t-}) dt,
\end{align*}
and their magnitude follows a conditional distribution $\mathrm{Gamma}(P_t,1)$,
\begin{align*}
    \mathbb{P} \big[ \zeta_n \in B \big\vert \sigma \{ (\tau_1, \zeta_1), \ldots, (\tau_{n-1}, \zeta_{n-1} ) \} \vee \sigma\{ \tau_n \} \big]
    &=
    \frac{ \Lambda( \{ \tau_n\}, B) }{ \Lambda( \{ \tau_n\},\mathbb{R}_{\geq 0} ) }
    \\
    &=
    \int_{ B } \frac{1}{\Gamma(P_{\tau_n})} \zeta^{P_{\tau_n} - 1} e^{-\zeta}d\zeta.
\end{align*}

The reason behind the specific $\Lambda$ in \eqref{Eq:Compensator-Example} comes from the use of the Gamma density in point processes for modelling natural disasters; see for example \cite{lesage_hawkes_2022}, where Hawkes processes with a Gamma density are considered for the modelling of insurance claims subjected to natural disasters. From these closed expressions planers can evaluate their position more accurately since now they know that, in this particular case, the  \textit{expected} destruction of any forthcoming catastrophe grows proportionally (in fact, linearly) to the state of pollution during the previous disaster.

\subsubsection{An intermediate model: PRMs with no diffusive term}
\label{Section:An_intermediate_model}

\n 

A quick inspection of this framework show that the current formulation indeed generalizes the previous models: the model in Section~\ref{section:BM+NHPP-driven-model} can be recovered as the special case in which the mark space is trivial and the intensity depends only on the pollution stock.

Another interesting case can be recovered by suppressing the diffusive term in $P$:
\begin{align}
    \label{Eq:K}  
	dK_t 
        &= 
        b^{\mathrm{cap}}(K_t,P_t,C_t,\theta_t)\,dt 
        - 
        \int_{ [0,\infty) } \big( 1 - \omega( K_{t-}, P_t, \zeta ) \big) K_{t-}\,q( dt, d\zeta),
    \\
    \label{Eq:P}
    dP_t 
        &= 
         b^{\mathrm{pol}}(K_t,P_t,C_t,\theta_t)\,dt 
    \\
    \label{Eq:K_0_P_0}
    K_0 & > 0,\quad P_0  > 0.
\end{align}
Intuitively, this model can be regarded as a natural extension to the one based on NHPP-driven model from Section \ref{section:Pollution-Driven-Disaster-Intensity-via-NHPP}. For the sake of presentation, we shall not focus on the interpretation of this model, as it can be inferred directly from the dynamics in \eqref{Eq:K-BM}--\eqref{Eq:K_0-P_0-BM}. However, we will refer to it as an intermediate step when deriving the HJB equation for the general jump-diffusion case.

\section{The Hamilton-Jacobi-Bellman PIDEs via the Dynamic Programming Principle}
\label{Section:Hamilton-Jacobi-Bellman-PIDEs}

\n 

\subsection{Preliminaries}

\n

Before proceeding, it is necessary to make a technical remark regarding the nature of the jumps under consideration, as well as the overall structure of the control problem. 

\subsubsection{Assumptions on the Stochastic Basis}
\label{Section:Assumptions_on_the_Basis}

\n 

As mentioned in Section \ref{section:Pollution-Driven-Disaster-Intensity-via-NHPP} and Remark \ref{Remark:Random-Intensity}, by considering a state-dependent intensity on the model, an implicit stochastic dependence is placed in the jumps in the form of a self-exciting interaction. This is possible because the jumps we are considering can all be derived from a Poisson point process, and more general, a Poisson random measure. 

To be more precise, the \textit{superposition} and \textit{thinning} properties of Poisson point processes allow us to assume, without any loss of generality,  that a marked process $q$ can be obtained as the integral of a larger counting random measure $N$:
\begin{align*}
    &q(dt \times d \zeta) =  \int \mathbf{1}_{ [0,\lambda_t(\zeta)] }(r) N( dt \times d \zeta \times dr \big),
\end{align*}
where $\lambda$ is a suitable (integrable, predictable) non-negative random field; see \cite{kingman_poisson_1993}, \cite{bremaud_point_2020} and \cite{daley_introduction_2008} for a reminder on general Poisson processes and random measures. Technically, this means we are working with jump processes of stochastic intensity. 

Whenever the source of randomness consists only on the jumps themselves, i.e. Sections \ref{section:Model under the constant arrival rate: standard Poisson Process}, \ref{section:Pollution-Driven-Disaster-Intensity-via-NHPP} and \ref{Section:An_intermediate_model}, our approach is justified by the embedding theorems on the extended state space $(K,P,\lambda)$ and the integration with respect to the Poisson random measure $N$ of intensity $dt \times dr$ on $\mathbb{R}_{\geq 0} \times \mathbb{R}_{\geq 0}$, see Chapter 5 of \cite{bremaud_point_2020}.

As we move to more complex specifications and introduce additional sources of randomness such as an independent Brownian motion, the stochastic basis $(\Omega, \mathcal{F}, \mathbb{F}, \mathbb{P} )$ must be correspondingly enlarged. Importantly, this will not pose a problem in what follows, since the Wiener-Poisson structure of the model, together with the additive nature of the proposed intensity $\lambda$, falls within the framework of \cite{hernandez-hernandez_coupled_2025} and \cite{hernandez-hernandez_mean-field_2026} for stochastic optimal control problems with environment-dependent jumps. Hence, the following assumption will be taken implicitly throughout the rest of the paper:

\begin{assumptions}
\label{Assumption-Poisson-Imbedding}
\sloppy
There exists an underlying Poisson random measure $N$ on $\mathbb{R}_{\geq 0} \times \mathbb{R}_{\geq 0} \times \mathbb{R}_{\geq 0}$ of intensity measure $dt \otimes \nu(d \zeta) \otimes dr$ for some $\sigma$-finite measure $\nu$ with finite second moment, such that
\begin{align}
    &q(dt,d \zeta) =  \int \mathbf{1}_{ [0,\lambda(P_t, \zeta)] }(r) N( dt, d\zeta, dr \big).
\end{align}
Additionally, we assume $\mathbb{F} = \{ \mathcal{F}_t \}_{t \geq 0}$ is (the complete, right-continuous augmentation of) the filtration
\begin{align}
    \mathcal{F}_t
    :=&
    \sigma \big\{ 
        q( B \times C \times D ) 
        ~\big|~ 
        (B,C,D) \in \mathcal{B}([0,t]) \otimes \mathcal{B}( \mathbb{R}_{\geq 0}) \otimes \mathcal{B}( \mathbb{R}_{\geq 0}) 
    \big\} \vee \mathcal{G}_t,
\end{align}
where
\begin{align*}
    \mathcal{G}_t
    =
    \sigma \big\{ 
        W_s
        ~\big|~ 
        0 \leq s \leq t
    \big\}
\end{align*}
if a Brownian motion is present in the model, and 
\begin{align*}
    \mathcal{G}_t
    =
    \{ \emptyset, \Omega \}
\end{align*}
otherwise.
\end{assumptions}

\subsubsection{Objective Functional}
\label{Sec:Objective-Functional}

\n

We now define the set of admissible control processes $C = (C_t)_{t\geq 0}$ and $\theta = (\theta_t)_{t \geq 0}$. To that end, we denote the set of \textit{admissible actions} as 
\begin{align} 
    \label{Eq:Admissible-actions}
    \mathfrak{a} := \mathbb{R}_{\geq 0} \times [0,\overline{\theta}],
    \quad
    \text{where}
    \quad
    \bar\theta:=\min\{1,\phi/z\}.
\end{align}
Particularly, observe that $E_t=\phi Y(K_t)-z\theta_t Y(K_t)\ge0$ for every $\theta \in [0,\overline{\theta}]$.

\begin{definition}[Admissible controls]
We call a control pair 
\begin{align*}
    (C_t,\theta_t)_{t \geq 0}=\big((C_t)_{t\ge0},(\theta_t)_{t\ge0}\big)
\end{align*}
\emph{admissible} for $(K,P) \in \mathcal{S}$ if
\begin{enumerate}

	\item 
        $(C_t,\theta_t)_{t \geq 0}$ is $\mathbb F$-predictable, with $(C_t,\theta_t) \in \mathfrak{a}$ $\mathbb{P}$-a.s. for every $t \geq 0$;

    \item 
        \sloppy
        under $(C_t,\theta_t)_{t \geq 0}$, the model's state equations admit a (pathwise) unique strong solution with nonnegative paths and initial condition $$(K_0,P_0) := (K,P);$$
	
    \item 
        the discounted utility is integrable:
    	\begin{equation}\label{eq:Admissible_Set}
    		\mathbb E \left[\int_0^\infty e^{-\rho t}\,\big|U(C_t,P_t)\big|\,dt\right]<\infty.
    	\end{equation}
\end{enumerate}
We denote the set of admissible controls by $\cA(K,P)$. 
\end{definition}

\begin{remark}
The expectation in \eqref{eq:Admissible_Set} is conditional on the initial state $(K,P)$, and $\cA(K,P) \neq \emptyset$ for all $(K,P) \in \cS$.
\end{remark}

For any admissible control $(C_t,\theta_t)_{t \geq 0} \in \cA(K,P)$, the associated gain function is defined by
\begin{equation}
    \label{Eq:Objective-Functional}
	J(K,P; C, \theta) := \bbE \left[ \int_{0}^{\infty} e^{-\rho s} U(C_s, P_s) \, ds \right],
\end{equation}
and the corresponding \textit{value function} $v: \cS \to \bbR$ is given by
\begin{equation} \label{eq:value_function}
	v(K,P) 
    := 
    \sup_{(C_t,\theta_t)_{t \geq 0} \in \cA(K,P)} J(K,P; C, \theta)
\end{equation}

In order to ensure the well-posedness of the problem and to derive the associated Hamilton-Jacobi-Bellman (HJB) equation in classical form, we impose the following conditions on the value function $v \colon \cS \to \mathbb{R}$:
\begin{enumerate}
    
    \item 
        \textit{Sufficient regularity:} $v\in C^2(\cS)$;
    
    \item 
        \textit{Monotonicity:}
        \begin{equation}
            \frac{\partial v}{\partial K} (K,P)\ge0,
            \qquad
            \frac{\partial v}{\partial P} (K,P)\le0,
            \qquad \forall (K,P)\in\cS;
        \end{equation}
    
    \item 
        \textit{Polynomial growth:} there exist constants $c>0$ and $\gamma,\delta\ge0$ such that
        \begin{equation}
            \label{eq:polynom_growth_v}
            |v(K,P)| \le c\big(1+K^\gamma+P^\delta\big),
            \qquad \forall (K,P)\in\cS.
        \end{equation}
\end{enumerate}

\subsubsection{Dynamic Programming Principle}

\n

Let $\cT_{t, T}$ denote the set of stopping times with values in the interval $[t, T]$, and define $\mathcal{T} := \cT_{0, \infty}$ as the set of admissible stopping times on the infinite horizon. The Dynamic Programming Principle (DPP) asserts that for any admissible initial condition $(K,P) \in \cS$, we have:
\begin{equation} 
    \label{eq:DPP}
	v(K,P) 
    = 
    \sup_{(C_t,\theta_t)_{t \geq 0} \in \cA(K,P)} 
        \sup_{h \in \cT}
        \bbE \left[ \int_0^h e^{-\rho s} U(C_s, P_s) \, ds + e^{-\rho h} v(K_h, P_h) \right].
\end{equation}
Equivalently, the value function also satisfies:
\begin{equation}
	v(K,P) 
    = 
    \sup_{(C_t,\theta_t)_{t \geq 0} \in \cA(K,P)} 
        \inf_{h \in \cT}
            \bbE \left[ \int_0^h e^{-\rho s} U(C_s, P_s) \, ds + e^{-\rho h} v(K_h, P_h) \right],
\end{equation}
with the convention that $e^{-\rho s(\omega)} = 0$ whenever $s(\omega) = \infty$.
This principle reflects the fact that optimal decision-making is time-consistent: the planner optimally balances immediate utility against the continuation value of the system's future state.

Suppose the candidate value function $v\in C^2(\cS)$ satisfies \eqref{eq:polynom_growth_v} and 
\begin{equation}
    \lim_{T\to\infty} \mathbb E  \left[e^{-\rho T} v(K_T,P_T)\right] 
    = 
    0.
\end{equation}
These conditions ensure that It\^o's formula applies and that the transversality condition holds.

Let $(C_t,\theta_t)_{t \geq 0} \in \cA(K,P)$ be admissible controls. Consider the system over a small time interval $[0, h]$, where $h > 0$. From the dynamic programming principle \eqref{eq:DPP}, we have the inequality:
\begin{equation} \label{eq:DPP_ineq}
	v(K,P) \geq \bbE \left[ \int_0^h e^{-\rho s} U(C_s, P_s) \, ds + e^{-\rho h} v(K_h, P_h) \right].
\end{equation}

\subsection{Analysis of model-specific equations}
\label{Section:Preliminaries-Solutions}

\n 

The Hamilton-Jacobi-Bellman (HJB) equation provides the infinitesimal version of the dynamic programming principle and characterizes the value function via a nonlinear partial integro-differential equation. We analyze the equation for each of the models presented earlier, under specific model-inspired assumptions. As stated in the Introduction, all corresponding computations are deferred to Appendix \ref{Sec:Technical-appendix} for the sake of readability.
Furthermore, the technical smoothness condition $v\in C^2(\cS)$ is assumed in order to work under the notion of so-called \textit{classical solutions}. We defer the discussion of viscosity solutions to Section~\ref{Sec:Verification-theorems-viscosity-solutions}.

\subsubsection{Unifying components}
\label{Section:Unifying-Components}

\n 

Since the framework encompasses models with similar characteristics, we follow the same roadmap as in the previous sections, starting from the simplest specification and mentioning only the changes that arise as the complexity increases. It is important to note that under the standing assumptions, the structure of the control problem cannot be altered substantially (in particular, the control of the agent is restricted to the trend or drift of the system and never on the noise or jump component). As a consequence, many of the results obtained for the base model extend naturally to the subsequent specifications. This is due to the fact that the Hamiltonian is preserved, which in turn maintains most of the structure of the optimality conditions for the controls. 

\subsubsection*{Hamiltonian and Optimality Conditions}
\label{Section:Hamiltonian and Optimality Conditions}

\n 

In each case, the control variables $(C, \theta)$ never directly acts outside the drift term in the dynamics of the state $(K,P)$, i.e. the jumps in capital and the diffusion in pollution are allowed to evolve uncontrolled. When translated to the HJB, this means that both the nonlocal and the second-order terms can be placed outside the supremum. This provides a cleaner representation of the dependence on the controls, since we can define the effective \textit{Hamiltonian} (i.e., only the part that depends on the controls) as
\begin{equation} 
\label{eq:Hamiltonian_StPois}
\begin{aligned}
	H(C,\theta;K,P,v_K,v_P) 
    :&= 
    U(C,P) 
        + v_K \, b^{\mathrm{cap}}( K, P, C, \theta ) 
        + v_P \, b^{\mathrm{pol}}( K, P, C, \theta ) 
\end{aligned}
\end{equation}
for all $v_K \geq 0$ and all $v_P \leq 0$. Furthermore, when referring to the \textit{maximized Hamiltonian} we shall use the notation
\begin{align}
    \label{eq:Hamiltonian_optimal}
    \widehat H(K,P)
    :=
    \sup_{(C,\theta) \in \mathfrak{a}} 
        H \Big( C,\theta; K,P, \frac{\partial v}{\partial K}(K,P), \frac{\partial v}{\partial P}(K,P) \Big).
\end{align}

Let $(\hat{C}_t,\hat{\theta}_t)_{ t \geq 0 }$ be the maximizers of $H$ over $\cA(K,P)$. On the one hand, observe that the optimal consumption $\hat C$ solves the interior first-order conditions:
\begin{equation} 
    \label{eq:HJB_FOC_C}
	\begin{aligned}
        \frac{\pt H}{\pt C}(\hat C, \hat \theta) 
        = \frac{\pt U}{\pt C}(\hat C, P) - \frac{\partial v}{\partial K}(K, P) = 0
        \iff 
        \frac{\pt U}{\pt C}(\hat C, P) = \frac{\partial v}{\partial K}(K, P).
    \end{aligned}
\end{equation}
On the other hand, from \eqref{eq:Z_I} observe that the optimal abatement share $\hat{\theta}$ satisfies the interior first-order condition
\begin{equation} 
    \label{eq:HJB_FOC_theta}
	\begin{aligned}
        \frac{\pt H}{\pt \theta} (\hat C, \hat \theta) 
        = 
        - \frac{\partial v}{\partial K}(K,P) Y(K) 
            - \frac{\partial v}{\partial P}(K,P)\, z\, Y(K) 
        = 
        0
        \\
        \iff 
        - \frac{\partial v}{\partial K}(K,P) = z \frac{\partial v}{\partial P}(K,P).
	\end{aligned}
\end{equation}

Consequently, the maximizer $\hat \theta$ is determined by the sign of \eqref{eq:HJB_FOC_theta} as follows:
\begin{equation} 
    \label{eq:Opt_theta}
	\hat \theta=
	\begin{cases}
		0,                                        & \tfrac{\partial v}{\partial K}{\scriptstyle (K,P)} + z \tfrac{\partial v}{\partial P}{\scriptstyle (K,P)} > 0,
        \\
		\text{any } \theta \in [0,\bar\theta],    & \tfrac{\partial v}{\partial K}{\scriptstyle (K,P)} + z \tfrac{\partial v}{\partial P}{\scriptstyle (K,P)} = 0,
        \\
		\bar\theta,                               & \tfrac{\partial v}{\partial K}{\scriptstyle (K,P)} + z \tfrac{\partial v}{\partial P}{\scriptstyle (K,P)} < 0,
	\end{cases}
\end{equation}

\subsubsection*{Candidate form for the value function}
\label{Section:Candidate form for the value function}

\n

In addition to the shared form of the Hamiltonian, we can further exploit our proposed framework in order to derive a closed-form expression for the value function $v$ which solves the dynamic programming equation~\eqref{eq:DPP}. 

To this end, we fix the following constants:
\begin{equation}
\psi \in \bbR_{>0}, \quad x \in \bbR_{>0}.
\end{equation}
Motivated by the first-order conditions \eqref{eq:HJB_FOC_C}--\eqref{eq:Opt_theta} and the dynamic programming
principle~\eqref{eq:DPP}, we postulate that for all $(K,P) \in \cS$ the partial derivatives of $v$ satisfy
\begin{align}
    \label{eq:value_fn_partial_K}
	\frac{\partial v}{\partial K}(K,P) = (\psi K)^{-\varepsilon}, 
    \quad
	\frac{\partial v}{\partial P}(K,P) = -x P^\beta.
\end{align}
From \eqref{eq:value_fn_partial_K}, integrating $\frac{\partial v}{\partial K}$ with respect to $K$ (for a fixed $P$) yields
\begin{align}\label{eq:v_from_K}
  v(K,P)
  = \psi^{-\varepsilon} \frac{K^{1-\varepsilon}}{1-\varepsilon}
  + c_1(P),
  \qquad \forall (K,P) \in \cS,
\end{align}
for some function $c_1 : \bbR_{>0} \to \bbR$. Similarly, from the integral $\frac{\partial v}{\partial P}$ with respect to $P$,
\begin{align}\label{eq:v_from_P}
  v(K,P)
  = -x \frac{P^{1+\beta}}{1+\beta}
  + c_2(K),
  \qquad \forall (K,P) \in \cS
\end{align}
for some function $c_2 : \bbR_{>0} \to \bbR$.

We introduce the functions $f,g : \cS \to \bbR$ by
\begin{align}
  f(K,P) := \psi^{-\varepsilon} \frac{K^{1-\varepsilon}}{1-\varepsilon}
  + c_1(P), \label{eq:f_partial_derivative}
  \quad
  g(K,P) := -x \frac{P^{1+\beta}}{1+\beta}
  + c_2(K). 
\end{align}
Therefore, the candidate solution $v$  is such that for all $(K,P) \in \cS$,
\begin{equation}
  v(K,P) = f(K,P) = g(K,P).
\end{equation}
Since $v$ is continuously differentiable on $\cS$, the same holds for $f$ and $g$, and their partial derivatives must coincide with those of $v$: for all $(K,P)\in \cS$, it must hold that
\begin{align} \label{eq:f_v_partial_P}
  \frac{\partial f}{\partial P}(K,P) = \frac{\partial v}{\partial P}(K,P),
  \quad
  \frac{\partial g}{\partial K}(K,P) = \frac{\partial v}{\partial K}(K,P).
\end{align}
We now use \eqref{eq:f_v_partial_P} together with \eqref{eq:f_partial_derivative} to determine the unknown functions $c_1$ and $c_2$. From \eqref{eq:v_from_K} we obtain
\begin{equation}
\frac{\partial f}{\partial P}(K,P)
= \frac{\partial}{\partial P}
\left(
\psi^{-\varepsilon} \frac{K^{1-\varepsilon}}{1-\varepsilon}
+ c_1(P)
\right)
= \frac{dc_1}{dP},
\end{equation}
for all $(K,P) \in \cS$.
Combining this with \eqref{eq:value_fn_partial_K} and \eqref{eq:f_v_partial_P} gives
\begin{align*}
    &\frac{dc_1}{dP} = \frac{\partial v}{\partial P}(K,P) = -x P^\beta
\end{align*}
for all $P> 0$. Thus $c_1$ satisfies the ordinary differential equation
\begin{equation*}
    \frac{dc_1}{dP} = -x P^\beta,
\end{equation*}
whose general $C^1$ solution is
\begin{equation}
    c_1(P) = -x \frac{P^{1+\beta}}{1+\beta} + C_1,
\end{equation}
for some constant $C_1 \in \bbR$. And similarly, for $\frac{\partial v}{\partial K}$ and $g$ in \eqref{eq:value_fn_partial_K} and \eqref{eq:f_v_partial_P},
\begin{equation}
c_2(K) = \psi^{-\varepsilon} \frac{K^{1-\varepsilon}}{1-\varepsilon}
+ C_2,
\end{equation}
for some constant $C_2 \in \bbR$.

Substituting $c_1$ and $c_2$ into~\eqref{eq:v_from_K}--\eqref{eq:v_from_P}, we find that for all $(K,P) \in \cS$,
\begin{equation*}
    v(K,P) 
    = 
    \psi^{-\varepsilon} \frac{K^{1-\varepsilon}}{1-\varepsilon} - x \frac{P^{1+\beta}}{1+\beta} + C_1
    = 
    \psi^{-\varepsilon} \frac{K^{1-\varepsilon}}{1-\varepsilon} - x \frac{P^{1+\beta}}{1+\beta} + C_2.
\end{equation*}
Hence, $v$ must be of the form
\begin{equation}
    v(K,P) 
    = 
    \psi^{-\varepsilon} \frac{K^{1-\varepsilon}}{1-\varepsilon} - x \frac{P^{1+\beta}}{1+\beta} + c,
\end{equation}
for some arbitrary constant $c = C_1 = C_2$. However, since the system \eqref{eq:value_fn_partial_K}, and thus the HJB equation and the associated optimal policies, are invariant under the addition of a constant to $v$, we may set $c = 0$ without loss of generality. As a result, the candidate value function $v$ reduces to
\begin{equation}\label{eq:value_fn}
  v(K,P)
  = \psi^{-\varepsilon} \frac{K^{1-\varepsilon}}{1-\varepsilon}
  - x \frac{P^{1+\beta}}{1+\beta},
  \qquad (K,P) \in \cS.
\end{equation}

Conversely, it is immediate to verify that the function $v$ defined in \eqref{eq:value_fn} satisfies \eqref{eq:value_fn_partial_K}. This justifies \eqref{eq:value_fn} as a natural candidate form for the value function. Moreover, taking into account \eqref{eq:Utility_Fn}, \eqref{eq:value_fn_partial_K}, \eqref{eq:value_fn}, as well as the first order conditions \eqref{eq:HJB_FOC_C} and \eqref{eq:HJB_FOC_theta}, we obtain
\begin{equation} 
    \label{eq:value_fn_derivatives_2}
	\begin{aligned}
		&\frac{\pt U}{\pt C}(\hat C, P) = \frac{\partial v}{\partial K}(K, P) 
        \iff 
        \hat{C}^{-\veps} = (\psi K)^{- \veps} 
        \iff 
        \hat{C} = \psi K.
	\end{aligned}
\end{equation}

\begin{remark}
In dynamic economic models, the equations in \eqref{eq:value_fn_partial_K} admit a natural interpretation as shadow prices: $\frac{\partial v}{\partial K}$ is the shadow value of an additional unit of capital and $\frac{\partial v}{\partial K}$ is the shadow cost of an additional unit of pollution. They are the dynamic-programming analogue of the adjoint equations from Pontryagin's Maximum Principle; see the well-known \cite{yong_stochastic_1999} for an in-depth discussion in the continuous-time stochastic case. 
\end{remark}

\subsubsection{Models of constant jump-size}

\n 

We begin with the base model introduced in Section \ref{section:Model under the constant arrival rate: standard Poisson Process}, where disasters are driven by a standard Poisson process of intensity $\lambda$. From the derivation presented in Appendix \ref{Sec:Technical-appendix}, the associated HJB equation takes the explicit form
\begin{equation} 
    \begin{aligned}
		\rho v(K, P) 
        = 
        \sup_{ (C,\theta) \in \mathfrak{a} } \bigg\{
            \frac{C^{1 - \varepsilon}}{1 - \varepsilon} 
            &- 
            \Big(\frac{\chi}{1 + \beta} - \alpha x \Big) P^{1 + \beta}
            + 
            \frac{(1 - \theta) AK - C}{ (\psi K)^{\varepsilon} }
            \\
            & - 
            x (\phi - \sigma \theta) A K P^\beta \
        \bigg\} +
        \lambda \big(v(\omega(K, P) K, P) - v(K, P) \big),&
	\end{aligned}
\end{equation}
for all $(K,P)\in \mathcal{S}$. Equivalently, using the notation for the optimal Hamiltonian from \eqref{eq:Hamiltonian_optimal}, we can write
\begin{align}
    \label{eq:HJB_max}
	\begin{aligned}
		\rho v(K,P) 
        =&~
        \widehat H(K,P)+\lambda\Big(v(\omega(K,P)K,P)-v(K,P)\Big).
	\end{aligned}
\end{align}
In other words, the nonlocality is driven by a \textit{homogeneous Poisson point process}.

Compare this with the next model, in which a \textit{nonhomogeneous Poisson} framework is used to model pollution-dependent disaster intensity: let $(\hat C, \hat \theta)$ denote the maximizers of $H$ over the admissible control set $\cA(K,P)$ as in \eqref{eq:Hamiltonian_optimal}. Since the jump component
\begin{equation} \label{eq:Nonhom_jump_term}
    \lambda(P)\big(v(\omega(K,P)K,P)-v(K,P)\big), 
\end{equation}
depends only on the state variables $(K,P)$, it is placed outside the Hamiltonian and does not affect the first-order conditions for the controls. That is, the HJB of the model from Section \ref{section:Pollution-Driven-Disaster-Intensity-via-NHPP} naturally decomposes into two components: 
\begin{enumerate}
\renewcommand{\theenumi}{\roman{enumi}}
    \item 
        the effective Hamiltonian term depending on the continuous dynamics and the controls, and 

    \item 
        the jump contribution, which depends on the state variables but not directly on the controls.
\end{enumerate}
In this case, the HJB equation takes the form
\begin{equation} \label{eq:HJB_max_NonHP}
	\rho v(K,P) = \widehat H(K,P) + \lambda(P)\Big(v(\omega(K,P)K,P)-v(K,P)\Big),
\end{equation}
for all $(K,P) \in \mathcal{S}$.

Moreover, given the fact that the first-order conditions \eqref{eq:HJB_FOC_C}--\eqref{eq:HJB_FOC_theta} characterize the feedback rules for optimal controls, by coupling them together with the pollution-dependent jump term in the HJB equation, see \eqref{eq:Nonhom_jump_term}, we obtain a system that determines the value function and the associated optimal policies. In relation to the homogeneous Poisson case (compare \eqref{eq:HJB_max_NonHP} with \eqref{eq:HJB_max}), the only structural change is the replacement of the constant intensity $\lambda$ by the state-dependent intensity $\lambda(P)=\lambda_0+\lambda_1 P$, see \eqref{eq:Nonhom_jump_term}. 

As discussed previously, this change introduces an additional endogenous feedback from pollution into the disaster hazard and, although deceptively simple, this modification carries deep technical implications. Fortunately, due to the choice of stochastic basis presented in Section \ref{Section:Assumptions_on_the_Basis}, we can naturally extend the previous setting in order to \textit{allow pollution to follow a diffusion process} while disaster arrivals remain governed by a nonhomogeneous Poisson process with intensity $\lambda(P_{t-})$. More precisely, following the product space construction in \cite{hernandez-hernandez_coupled_2025} and \cite{hernandez-hernandez_conditional_2025} -- see also the details in the case of control problems in \cite{hernandez-hernandez_mean-field_2026} -- as well as a similar derivation procedure as in the previous cases (see Appendix \ref{Sec:Technical-appendix} for the details), we obtain that the HJB equation for the model described in Section \ref{section:BM+NHPP-driven-model} is
\begin{equation} 
    \label{eq:HJB_max_NonHP+Diff}
	\rho v(K,P) 
    = 
    \widehat H(K,P) 
    +
    \frac{1}{2} \,\sigma^2\,P^2\,\frac{\partial^2 v}{\partial P^2}(K,P)
    + 
    \lambda(P)\Big(v(\omega(K,P)K,P)-v(K,P)\Big),
\end{equation}
for all $(K,P) \in \mathcal{S}$. 

Equation \eqref{eq:HJB_max_NonHP+Diff} above shows that the effect of disasters enters additively via the compensator-adjusted jump term. Thus the marginal damage of emissions is amplified through both the continuous deterioration of environmental quality (via the drift and diffusion of $P_t$) and the increased likelihood of discrete catastrophic events.

On the other hand, the second-order term 
\begin{align}
    \label{Eq:2nd-order-HJB-term}
    \frac{1}{2} \,\sigma^2\,P^2\,\frac{\partial^2 v}{\partial P^2}(K,P)
\end{align}
captures the effect of Brownian pollution shocks. The combination of diffusion and jump risk implies that optimal policies $(C,\theta)$ must balance the trade-off between consumption, abatement, and the endogenous exposure to both continuous and discontinuous environmental risks.

\subsubsection{Models of varying jump-size}

\n 

Moving on to the generalized framework of Poisson random measures -- recall the specifications made in Section \ref{Eq:PRM-Model-Introduction} -- the jump component of the state dynamics is no longer characterized by a single nonlocal term, but by an integral over all possible jump magnitudes. 

Proceeding as before, and using the compensated PRM in the dynamic programming argument -- see \eqref{Eq:Poisson-Measure-HJB} in the Appendix \ref{Sec:Technical-appendix} -- the HJB equation for the intermediate model from Section \ref{Section:An_intermediate_model} becomes
\begin{equation}
    \label{eq:HJB_max_PRM}
    \begin{aligned}
        \rho v(K,P) 
        =
        \widehat{H}(K,P) 
        +
        \int_{(0,\infty)} 
            \big( v( \omega( K,P, \zeta )K,P) - v(K,P) \big) \lambda(P,\zeta)
        \,\nu( d\zeta ),
    \end{aligned}
\end{equation}
for all $(K,P) \in \mathcal{S}$.

It is clear that equation \eqref{eq:HJB_max_PRM} above represents a direct generalization from nonhomogeneous Poisson model with deterministic jump size from \eqref{eq:HJB_max_NonHP}. In this formulation, \textit{disaster arrivals are described by a Poisson random measure with marks} $\zeta \in (0,\infty)$, where $\nu(\mathrm d\zeta)$ governs the distribution of jump characteristics and $\lambda(P,\zeta)$ allows the arrival intensity to depend jointly on pollution and on the mark itself. Relative to \eqref{eq:HJB_max_NonHP}, the jump contribution is no longer summarized by a single nonlocal term evaluated at a fixed loss factor. Instead, the expected impact of disasters is obtained by integrating the value loss
\begin{equation*}
    v\big(\omega(K,P,\zeta)K,P\big)-v(K,P)
\end{equation*}
over the entire mark space $(0,\infty)$. This structure captures heterogeneity in disaster severity while preserving the separation between the control-dependent Hamiltonian $\widehat H(K,P)$ and the purely state-dependent jump contribution. As in the previous cases, the first-order conditions characterizing the optimal controls are unaffected by the jump term, since neither $\omega(K,P,\zeta)$ nor $\lambda(P,\zeta)$ depends on $(C,\theta)$.

From an economic perspective, \eqref{eq:HJB_max_PRM} already embeds a rich feedback from pollution into disaster risk. Pollution may increase the frequency of arrivals, distort the distribution of marks toward more destructive events, or both. The value function therefore internalizes not only the expected number of disasters but also their expected severity, aggregated through the integral operator.

Lastly, due to our current formulation, the transition from \eqref{eq:HJB_max_PRM} to \textit{the general jump-diffusion setting} is both conceptually and technically straightforward: allowing pollution to follow a diffusion process adds the local second-order term from \eqref{Eq:2nd-order-HJB-term}, which captures continuous environmental risk. More precisely, 
\begin{equation}
    \label{eq:HJB_max_PRM+diff}
    \begin{aligned}
        \rho v(K,P) 
        =
        \widehat{H}(K,P) 
        &+ 
        \frac{1}{2} \,\sigma^2\,P^2\,\frac{\partial^2 v}{\partial P^2}(K,P)
        \hspace{10em}
        \\
        &+
        \int_{(0,\infty)} 
            \big( v( \omega( K,P, \zeta )K,P) - v( K,P ) \big) \lambda(P,\zeta) \nu( d\zeta ),
    \end{aligned}
\end{equation}
for all $(K,P) \in \mathcal{S}$, see \eqref{Eq:Brownian-HJB}.

Importantly, this extension leaves the jump operator unchanged: discrete disaster risk and continuous pollution risk enter additively in the HJB equation. The resulting value function reflects joint exposure to smooth fluctuations in environmental quality and to rare but potentially severe catastrophic shocks, both endogenously linked to the pollution state.

\section{Verification theorems and viscosity solutions}
\label{Sec:Verification-theorems-viscosity-solutions}

\n

In this section, we consider the \emph{verification problem}, that is, the conditions under which the solution of a Hamilton-Jacobi-Bellman partial integro-differential equation (PIDE) coincides with the value function of the underlying control problem.  Since it is well-known in the literature that the value function of an optimal control problem need not be smooth, we distinguish between two solution concepts: \emph{classical} and \emph{viscosity} solutions; see \cite{crandall_users_1992}. A classical solution is a sufficiently differentiable function that satisfies the PIDE pointwise.
In contrast, a viscosity solution is a function for which the notion of differentiability itself is weakened, and interprets the PIDE in a generalized sense.
Both notions are relevant in our setting, as the nonlocal jump terms in the HJB equations affect the regularity of possible solutions, especially as the complexity of the jump processes grow.

We formulate and prove the verification theorems for the previously introduced jump-control growth-environment models, based on the results from \cite{garroni_second_2002} on second-order integro-differential problems, and \cite{gichman_stochastic_1972} on general SDEs with jumps; for an in-depth treatment on the Brownian case we refer to \cite{pham_continuous-time_2009}.

\subsection{Verification theorem}

\n 

The central step in the classical dynamic programming approach is to establish that, whenever a sufficiently smooth solution to the Hamilton-Jacobi-Bellman equation exists, this solution coincides with the value function of the control problem. This result, known as the verification theorem, further yields the existence of an optimal Markovian control as a direct corollary. The proof relies fundamentally on the application of It\^o's formula.

Recall that the planner maximizes the objective functional $J(K, P ;C,\theta)$ defined in \eqref{Eq:Objective-Functional}, with the corresponding value function $v(K,P)$ as in \eqref{eq:value_function}, and let $C_0(\mathcal{S})$ denote the space of continuous functions vanishing at infinity with the supremum norm. From the theory of Markov processes, it is well known that for a given control $(C_t, \theta_t)_{t \geq 0}$ there exists a (pseudo-)differential operator 
\begin{align*}
    \mathcal{L}^{C,\theta} 
    : 
    \mathrm{Dom}(\mathcal{L}^{C,\theta}) \subset C_0(\mathcal{S}) 
    \to
    C_0(\mathcal{S})
\end{align*}
acting as the \textit{(controlled) infinitesimal generator} of $(K_t,P_t)_{t \geq 0}$, see \cite{fleming_controlled_2006}.

Furthermore, the HJB equation of $v$ can be written in terms of the generator as
\begin{equation} 
    \label{eq:HJB_VerTh} 
    \rho v(K,P)
    =
    \sup_{(C,\theta)\in \mathfrak{a}}\Big\{
        U(C,P)
        + 
        \mathcal L^{C,\theta}v(K,P)
    \Big\}.
\end{equation}

\sloppy
With this in mind, we employ the present framework to provide explicit verification results for HJBs in the form of \eqref{eq:HJB_VerTh} for each of the models under consideration. In particular, we prove that the value function $v$ constitutes an upper bound for the objective functional defined in Section \ref{Sec:Objective-Functional}, 
\begin{equation}
    \label{Eq:v-upper-bound}
    v(K,P)\ \ge\ J(K, P; C, \theta),
\end{equation}
with $J$ as in \eqref{Eq:Objective-Functional}; furthermore, we prove the existence of admissible controls $(\hat C_t, \hat \theta_t)_{t \geq 0} \in \mathcal{A}(K,P)$ that are in fact \textit{optimal}; that is, they satisfy the equality
\begin{equation}
    \label{Eq:Optimality-notion}
    v(K,P)
    =
    J(K,P; \hat C, \hat \theta)
    =
    \sup_{(C,\theta)\in\cA(K, P)}J(K,P;C,\theta).
\end{equation}

\begin{remark}
A key characteristic of the problem in an infinite-time horizon is that, unlike its finite counterpart, there is no specific terminal condition. In its place, we work with the so-called \textit{transversality condition},
\begin{equation}
    \label{eq:TV}
	\lim_{T\to\infty}\mathbb E \left[e^{-\rho T}v(K_T,P_T)\right] = 0,
\end{equation}
which, as we will see below, will be helpful when obtaining the optimality of the controlled processes.
\end{remark}

\subsubsection{State-independent Jumps (Baseline HPP-driven model)}

\n 

Observe that in the baseline model, the controlled generator is defined as
\begin{equation}
\label{eq:Gen-verify}
\begin{aligned}
	\mathcal L^{C,\theta}\varphi(K,P)
	:= 
    \frac{\partial \varphi}{\partial K}(K,P)\big((1-\theta)AK-C\big) 
    &+ 
    \frac{\partial \varphi}{\partial P}(K,P)\big((\phi-z\theta)AK-\alpha P\big)
    \\[1.5pt]
	&+ \lambda\big( \varphi(\omega(K, P)K,P) - \varphi(K,P)\big),
\end{aligned}
\end{equation}
for any test function $\varphi \in C^2(\cS)$.

\begin{theorem}[Verification: Jump-control growth-environment model]
\label{thm:Pham-verify} 
Suppose $v\in C^2(\cS)$ is nonnegative, has at most polynomial growth, and satisfies the HJB equation \eqref{eq:HJB_VerTh} pointwise, with $\mathcal L^{C,\theta}\varphi(K,P)$ as in \eqref{eq:Gen-verify}.

\smallskip
\item[\emph{(i)}]
    \sloppy
    Upper bound. 
    Let $(K_t,P_t)_{t \geq 0}$ be the state process of an economy driven by \eqref{Opt_K}--\eqref{Opt_P} for some admissible control $(C_t, \theta_t)_{t \geq 0}\in\mathcal A(K, P)$. If the transversality condition \eqref{eq:TV} holds, then equation \eqref{Eq:v-upper-bound} holds.

\smallskip
\item[\emph{(ii)}]
    Optimality. If there exist Borel functions $(\hat C,\hat \theta) : \mathcal{S} \to \mathfrak{a}$ such that
	\begin{equation} 
        \label{eq:Att}
		U(\hat C,P)+ \mathcal L^{\hat C(K,P), \hat\theta(K,P)}v(K,P)
		=
        \sup_{(C,\theta)\in \mathfrak{a} }\Big\{U(C,P) + \mathcal L^{C,\theta}v (K,P)\Big\},
	\end{equation}
    for all $(K,P) \in \mathcal{S}$ with $\mathcal L^{C,\theta}\varphi(K,P)$ as in \eqref{eq:Gen-verify}, then the admissible control 
    \begin{align}
        \label{Eq:Generic-Feedback-Control}
        (\hat C_t, \hat \theta_t)_{t \geq 0} 
        := 
        \big( \hat C(K_t,P_t),\hat \theta(K_t,P_t) \big)_{t \geq 0} 
        \in 
        \mathcal{A}(K,P)
    \end{align}
    is \textit{optimal}, in the sense of \eqref{Eq:Optimality-notion}.
\end{theorem}

\begin{proof}

\item[(i)] 
    Fix any admissible $(C_t,\theta_t)$ with state $(K_t, P_t)$. Define $Y_t := e^{-\rho t}v(K_t, P_t)$. Applying It\^o's formula for jump processes to $Y_t$ on $[0, T]$, where $T < 0$ is an arbitrary time horizon, gives
	\begin{equation}
	\begin{aligned}
		e^{-\rho T}v(K_T,P_T)-v(K_0,P_0)
		&=\int_0^T e^{-\rho s}\big(\cL^{C_s,\theta_s} v(K_s,P_s)-\rho v(K_s,P_s)\big)\,ds + M_T,
	\end{aligned}
	\end{equation}
    where 
    \begin{equation} 
        \label{eq:M_T}
        M_T 
        := 
        \int_0^T e^{-\rho s}\Big( 
            v(\omega(P_{s-},K_{s-})K_{s-},P_{s-}) - v(K_{s-},P_{s-})
        \Big)\,d\tilde q_s
    \end{equation}
    is a martingale. Equivalently, \eqref{eq:M_T} can be written as
    \begin{equation}
    \begin{aligned}
    	M_T
    	&=\sum_{0<\tau_n\le T} e^{-\rho \tau_n}\Big(v(\omega(P_{\tau_n-},K_{\tau_n-})K_{\tau_n-},P_{\tau_n-})-v(K_{\tau_n-},P_{\tau_n-})\Big)\\
    	&\quad -\lambda\int_0^T e^{-\rho s}\Big(v(\omega(P_s, K_s)K_s,P_s)-v(K_s,  P_s)\Big)\,ds.
    \end{aligned}
    \end{equation}
    Taking expectations and using $\bbE[M_T] = 0$,
    \begin{equation}
        \label{Eq:Feynam-Kac}
        \mathbb E \left[e^{-\rho T}v(K_T,P_T)\right]-v(K_0, P_0)
        =
        \mathbb E \left[\int_0^T e^{-\rho s}\big( 
            \cL^{C_s,\theta_s}v(K_s,P_s) - \rho v(K_s,P_s)
        \big)\,ds\right].
    \end{equation}
    By \eqref{eq:HJB_VerTh}, we obtain
    \begin{equation}
    \rho v \ge U(C_t,P_t)+\mathcal L^{C_t,\theta_t}v \iff \mathcal L^{C_t,\theta_t}v-\rho v \le -U(C_t,P_t).
    \end{equation}
For every $T$,
    \begin{equation}
    \mathbb E \left[e^{-\rho T}v(K_T,P_T)\right]-v(K_0,P_0)
    \ \le\ -\,\mathbb E \left[\int_0^T e^{-\rho t}U(C, P)\,dt\right].
    \end{equation}
    
    Then, through a localization argument (see Section 7.2 of \cite{garroni_second_2002}, Theorem 7.2.1, for the solution of the discounted control problem for $(t, (K,P))$ on the localized domain $[0,T] \times [-n, n]^2$ for any $n \geq 1$), we can use standard results on the asymptotic stability of Wiener-Poisson semigroups to extend the (sub-optimally) controlled generator from \eqref{eq:Gen-verify} up to $[0,\infty) \times \mathcal{S}$ (see Theorem 1 of III.13 from \cite{gichman_stochastic_1972}); as a result, we can let $T\to\infty$ so that \eqref{eq:TV} becomes 
    \begin{equation}
    -v(K,P)\ \le\ -J(K,P;C,\theta),\qquad\text{i.e.}\quad v(K,P)\ \ge\ J(K,P;C,\theta).
    \end{equation}   
\item[(ii)]	
Let $(\hat C, \hat \theta)$ satisfy \eqref{eq:Att}. Repeating the above with $(\hat C, \hat \theta)$, the HJB inequality becomes an equality:
\begin{equation}
\rho v(K,P) = U(\hat C,P)+\mathcal L^{\hat C, \hat \theta}v(K,P).
\end{equation}
So for each $T$,
\begin{equation}
\mathbb E \left[e^{-\rho T}v(\hat K_T, \hat P_T)\right]-v(K,P)
= -\,\mathbb E \left[\int_0^T e^{-\rho t}U(\hat C_t, \hat P_t)\,dt\right].
\end{equation}
As before, we use a localization argument to let $T\to\infty$ (this time for the optimally controlled generator); combining this with \eqref{eq:TV} for $(\hat C, \hat \theta)$ yields $v(K,P)=J(K,P; \hat C, \hat \theta)$. Combined with (i) this implies optimality.
\end{proof}

\begin{corollary}[Verification for the candidate power-separable value function]
Let $v$ be given by
\begin{equation}
    \tag{\ref{eq:value_fn}}
	v(K,P)=\psi^{-\varepsilon}\frac{K^{1-\varepsilon}}{1-\varepsilon}-x\frac{P^{1+\beta}}{1+\beta},
\end{equation}
with $x > 0$ and $\psi>0$. Let $\hat C$ and $\hat \theta$ be the feedback controls defined by \eqref{eq:value_fn_derivatives_2} and \eqref{eq:Opt_theta}, respectively. Then for all $(K,P)\in \mathcal{S} $,
\begin{equation}
	v(K,P) = \sup_{(C,\theta)\in\mathcal A(K,P)}J(K,P;C,\theta) = J(K,P;\hat C, \hat \theta).
\end{equation}
\end{corollary}

\begin{proof}
{
We divide the argument into several steps in order to verify  that the candidate function $v$ is indeed the value function of the control problem:

\medskip
\item[1.] 
    \textit{Regularity and growth of $v$.} Let $\varepsilon,\beta>0$ so that $v\in C^2(\mathcal{S})$. Then, in addition to the first-order conditions from \eqref{eq:value_fn_partial_K}, it also holds that
    \begin{equation}
    \begin{aligned}
        \frac{\partial^2 v}{\partial K^2}(K,P) = -\varepsilon \psi^{-\varepsilon}K^{-\varepsilon-1},\,
        \frac{\partial^2 v}{\partial P^2}(K,P) = -x\beta P^{\beta-1},\,
        \frac{\partial^2 v}{\partial K \partial P}(K,P) = 0.
    \end{aligned}
    \end{equation}
    Hence $v$ has at most polynomial growth in $(K,P)$ and it satisfies the regularity requirements from Theorem~\ref{thm:Pham-verify}.

\medskip
\item[2.] 
    Although the feedback rule $\hat\theta$ is of \textit{bang-bang type} and introduces a discontinuity in the drift across the switching surface $\{(K,P):\frac{\partial v}{\partial K}+z \tfrac{\partial v}{\partial P}=0\}$, existence and uniqueness of strong solutions remain valid. Indeed, the drift is piecewise affine and thus of bounded variation on compacts. Consequently, using Theorem 4.1 from \cite{przybylowicz_existence_2021} for the well-posedness of SDEs with discontinuous drifts, $(K_t,P_t)$ has a unique global strong existence. 

\medskip
\item[3.] 
    \textit{Maximization of the Hamiltonian.} On the one hand, for fixed $(K,P)$ the map $C\mapsto H$ is strictly concave, with $\hat C=\psi K$ as the first-order condition for optimality, see \eqref{eq:value_fn_derivatives_2}. On the other hand, since the map $\theta\mapsto H$ is affine with slope $-AK\,(\frac{\partial v}{\partial K}+z \tfrac{\partial v}{\partial P})$, its supremum on $[0,\bar{\theta}]$ is attained at $\hat{\theta}$, defined in \eqref{eq:Opt_theta}. Hence, for each $(K,P) \in \mathcal{S}$,
    \begin{equation}\label{eq:attainment-pointwise}
    	\begin{aligned}
    		\sup_{ (C, \theta) \in \mathfrak{a} } \left\{ 
                \frac{C^{1-\veps}}{1-\veps} 
                -
                \chi\frac{P^{1+\beta}}{1+\beta} 
                + 
                \mathcal L^{C,\theta}v(K,P)
            \right\} 
            =
            \frac{\hat C^{1-\veps}}{1-\veps} 
            - 
            \chi\frac{P^{1+\beta}}{1+\beta} 
            +
            \mathcal L^{\hat C, \hat \theta}v(K,P);
    	\end{aligned}
    \end{equation}
    i.e. $(\hat C,\hat\theta)$ is a maximizer of the Hamiltonian \eqref{eq:Hamiltonian_StPois}.

\medskip
\item[4.] 
    \textit{Transversality.} Since $v$ has polynomial growth
    \begin{equation}
    v(K,P)\le c_1(1+K^{1-\varepsilon}+P^{1+\beta})
    \end{equation}
    for some $c_1 > 0$ and $(K_t,P_t)$ has moments with at most exponential growth rate strictly smaller than $\rho$ (this holds under linear drift/jump structure and $\theta\in[0, \bar{\theta}]$), we obtain that the condition
    \begin{equation}
    	\lim_{T\to\infty}\mathbb E \big[ 
            e^{-\rho T}v(K_T,P_T) \big| (K_0,P_0) = (K,P) 
        \big] 
        = 0
    \end{equation}
    holds for any initial state $(K,P) \in \cS$.

\medskip
\item[5.] 
    \textit{Upper bound via verification.} Let $(C_t,\theta_t)_{t \geq 0} \in \mathcal{A}(K,P)$. Observe that by applying It\^o's formula to $e^{-\rho t}v(K_t,P_t)$ and taking expectation, the transversality condition \eqref{eq:TV} and the upper bound \eqref{Eq:v-upper-bound} from Theorem \ref{thm:Pham-verify}.(i) yield
    \begin{equation}
    		v(K,P)\;\ge\;J(K,P;C,\theta)
    \end{equation}
    Then, since $(\hat C, \hat \theta)$ attains the supremum of the Hamiltonian by definition, the previous inequality becomes an equality.

    Combining this with Step 4 gives the stated identity. Thus, for every initial state $(K,P)\in\mathcal{S}$,
    \begin{equation}
    	v(K,P)
    	=
        \sup_{(C_t,\theta_t)_{t \geq 0} \in \cA(K,P)} \mathbb E \bigg[
            \int_0^\infty e^{-\rho t}\Big( 
                \frac{C_t^{\,1-\varepsilon}}{1-\varepsilon}
                -
                \chi\frac{P_t^{\,1+\beta}}{1+\beta}
            \Big) dt
        \bigg]
    	= 
        J\big(K,P;\hat C,\hat\theta\big).
    \end{equation}
}
\end{proof}

With the previous result \textit{we have established an explicit candidate for a power-separable value function $v$ and the corresponding optimal control policy $(\hat C,\hat\theta)$}, both satisfying all the conditions of the verification theorem.

\subsubsection{State-dependent Jumps (NHPP-driven model)}

\n 

We now state and prove a verification theorem for the planner's problem when the arrival of destructive events is governed by a NHPP with state-dependent intensity \eqref{eq:NonHP_lambda}. In this case, the controlled infinitesimal generator is given by
\begin{equation}
\label{eq:generator-nonHP}
\begin{aligned}
	\mathcal L^{C,\theta}\varphi(K,P)
	:= 
    \frac{\partial \varphi}{\partial K}(K,P)\big((1-\theta)AK-C\big) 
    &+ 
    \frac{\partial \varphi}{\partial P}(K,P)\big((\phi-z\theta)AK-\alpha P\big)
	\\[1.5pt]
	&+\lambda(P)\Big( \varphi(\omega(P,K)K,P) - \varphi(K,P)\Big)
\end{aligned}
\end{equation}
for any test function $\varphi\in C^2(\cS)$.

\begin{theorem}[Verification: jump-control growth-environment model with state-dependent intensity]
\label{thm:verification-nonHP}
Suppose $v\in C^2(\cS)$ is nonnegative, has at most polynomial growth, and satisfies the HJB equation \eqref{eq:HJB_VerTh} pointwise, with $\mathcal L^{C,\theta}\varphi(K,P)$ as in \eqref{eq:generator-nonHP}.

\smallskip
\item[\emph{(i)}]
    \sloppy
    Upper bound.
    Let $(K_t,P_t)_{t \geq 0}$ be the state process of an economy driven by \eqref{eq:Opt_K_NHPP}--\eqref{eq:Opt_P_NHPP} for some admissible control $(C_t, \theta_t)_{t \geq 0}\in\mathcal A(K, P)$. If the transversality condition \eqref{eq:TV} holds, then equation \eqref{Eq:v-upper-bound} holds.
		
\smallskip
\item[\emph{(ii)}]
    Optimality. 
    If there exist Borel functions $(\hat C,\hat \theta) : \mathcal{S} \to \mathfrak{a}$ such that \eqref{eq:Att} holds for all $(K,P) \in \mathcal{S}$ with $\mathcal L^{C,\theta}\varphi(K,P)$ as in \eqref{eq:generator-nonHP}, then the admissible control $(\hat C_t, \hat \theta_t)_{t \geq 0}$, constructed as in \eqref{Eq:Generic-Feedback-Control}, is \textit{optimal} in the sense of \eqref{Eq:Optimality-notion}.
    
\end{theorem}

\begin{proof}

Due to the Poissonian nature of the jumps, the result is obtained by following the proof of Theorem \ref{thm:Pham-verify} \textit{verbatim} except for a substitution of the jump martingale term by
\begin{equation}
    M_T 
    := 
    \int_0^T e^{-\rho s}\Big(v(\omega(K_{s-}, P_{s-})K_{s-},P_{s-})-v(K_{s-},P_{s-})\Big)\,d\tilde q_s.
\end{equation}

\end{proof}

\subsubsection{State-dependent Jump-Diffusions (NHPP-driven model)}

\n

In the present jump-diffusion setting with pollution-dependent intensity, the Hamiltonian retains the same $(C,\theta)$ dependence as in the constant-intensity case -- hence, the first-order conditions deliver the feedback rules $\hat C=\psi K$ and the $\theta$-projection in \eqref{eq:Opt_theta}-- with the main difference being that generator involved is now of second-order:
\begin{equation}
    \label{eq:Gen-verify-NH}
	\begin{aligned}
	\mathcal L^{C,\theta}\varphi(K,P)
	:=
    \frac{\partial \varphi}{\partial K}(K,P)\big((1-\theta)AK-C\big) 
    + 
    \frac{\partial \varphi}{\partial P}(K,P)\big((\phi-z\theta)AK-\alpha P\big)
	\\
	+
    \frac{1}{2} \sigma^2 P^2 \frac{\partial^2 \varphi}{\partial P^2}(K,P)
    +
    \lambda(P)\Big( \varphi(\omega(P,K)K,P) - \varphi(K,P)\Big),
    \end{aligned}
\end{equation}
for any test function $\varphi\in C^2(\cS)$.

However, as stated earlier, this technicality is surpassed due to our assumptions on the stochastic basis (Section \ref{Section:Assumptions_on_the_Basis}). With these specifications in place, Theorem~\ref{thm:verify-NH} below provides a direct generalization of Theorems \ref{thm:Pham-verify} and \ref{thm:verification-nonHP}.

\begin{theorem}[Verification: nonhomogeneous Poisson with Brownian pollution]
\label{thm:verify-NH}
Suppose $v\in C^2(\cS)$ is nonnegative, has at most polynomial growth, and satisfies the HJB equation \eqref{eq:HJB_VerTh} pointwise, with $\mathcal L^{C,\theta}\varphi(K,P)$ as in \eqref{eq:Gen-verify-NH}.

\smallskip
\item[\emph{(i)}]
    \sloppy
    \emph{Upper bound.} 
    Let $(K_t,P_t)_{t \geq 0}$ be the state process of an economy driven by \eqref{eq:Opt_K_NHPP} and \eqref{eq:P_W_NonH_Pois} for some admissible control $(C_t, \theta_t)_{t \geq 0}\in\mathcal A(K, P)$. If the transversality condition \eqref{eq:TV} holds, then equation \eqref{Eq:v-upper-bound} holds.
		
\smallskip
\item[\emph{(ii)}]
    \emph{Optimality.} 
    If there exist Borel functions $(\hat C,\hat \theta) : \mathcal{S} \to \mathfrak{a}$ such that \eqref{eq:Att} holds for all $(K,P) \in \mathcal{S}$ with $\mathcal L^{C,\theta}\varphi(K,P)$ as in \eqref{eq:Gen-verify-NH}, then the admissible control $(\hat C_t, \hat \theta_t)_{t \geq 0}$, constructed as in \eqref{Eq:Generic-Feedback-Control}, is \textit{optimal} in the sense of \eqref{Eq:Optimality-notion}.

\end{theorem}

\begin{proof}
Let $(C_t,\theta_t)$ admissible be fixed with the corresponding state $(K_t,P_t)$, and define $Y_t:=e^{-\rho t}v(K_t,P_t)$. By It\^o's formula,
\begin{equation}\label{eq:Ito-Y-NH}
	\begin{aligned}
		e^{-\rho T}v(K_T,P_T)-v(K_0,P_0)
		=& 
        \int_0^T e^{-\rho s}\Big((\mathcal L^{C_s,\theta_s}v)(K_s,P_s)-\rho v(K_s,P_s)\Big)\,ds 
            \\
			&+ 
            \underbrace{\int_0^T e^{-\rho s} \frac{\partial v}{\partial P}(K_s,P_s) \sigma_s\,dW_s}_{=:M_T^{(W)}} 
            + 
            \underbrace{\int_0^T e^{-\rho s} \Delta v_s\,d\tilde{q}_s}_{=:M_T^{(\hat q)}},
		\end{aligned}
	\end{equation}
where $\tilde{q}$ is the compensated jump martingale of $q$,
\begin{equation}
	\Delta v_s:=v(\omega(K_{s-},P_s)K_{s-},P_s)-v(K_{s-},P_s),
\end{equation}
and we have used
	\begin{equation}
	\int_0^T e^{-\rho s}\Delta v_s\,d q_s
	=
    \int_0^T e^{-\rho s}\Delta v_s\,d\tilde{q}_s
	   +
       \int_0^T e^{-\rho s}\lambda(P_s)\Delta v_s\,ds.
\end{equation}
Then, as in the proof of Theorem \ref{thm:verification-nonHP}, the result is obtained by following the proof of Theorem \ref{thm:Pham-verify} \textit{verbatim}, substituting $M$ with $M^{(W)} + M^{(\tilde{q})}$ from equation \eqref{Eq:Feynam-Kac} and onward.
\end{proof}

In terms of power-separable value functions, it remains to check that our candidate $v$ from \eqref{eq:value_fn} indeed solves the HJB and that $(\hat C,\hat\theta)$ attains the supremum. Substituting the candidate function into \eqref{eq:HJB_VerTh} with respect to the generator \eqref{eq:Gen-verify-NH}, the contributions of $K$ and $P$ can be treated separately. The resulting $v$ is $C^2$ with polynomial growth, and the induced feedback is admissible and attains the Hamiltonian pointwise. Hence, by Theorem~\ref{thm:verify-NH}, the candidate is the value function and the feedback is optimal.

\subsubsection{On the transversality condition}

\n 

As we have seen, the transversality condition \eqref{eq:TV} turned out to be crucial for the verification of solutions to the HJBs in each model. In practice, it may not be possible to check \eqref{eq:TV} directly; instead, one can derive a set of sufficient conditions for the transversality condition to hold. 

\begin{assumptions}\noindent

    \item[(i)] 
        \textit{Polynomial growth.} The value function $v$ has at most polynomial growth. 
        
        \smallskip
        \noindent
        - In addition, for the jump-diffusion case (Theorem \ref{thm:verify-NH}), it holds that $v\in C^2(\cS)$ with derivatives of at most polynomial growth: there exists a constant $C > 0$ such that for any $(K,P) \in \mathcal{S}$,
        \begin{equation} 
			\left|v(K,P)\right|
            +
            \left| \frac{\partial v}{ \partial K }(K,P) \right|
            +
            \left| \frac{\partial v}{ \partial P }(K,P) \right|
            +
            \left| \frac{\partial^2 v}{ \partial K^2 }(K,P) \right| 
            \le 
            C(1+K^m+P^m).
        \end{equation}
    
    \item[(ii)]
        \textit{Integrability of the state.} For any admissible control, there exists $m>0$ such that 
		\begin{equation}
			\sup_{t \geq 0} \mathbb{E}[|K_t|^m + |P_t|^m] < \infty.
		\end{equation}

    \item[(iii)]
        \textit{Integrability of the intensity.} For every $T > 0$, 
        \begin{align*}
            \bbE \Big[\int_0^T \lambda(P_s)\,ds\Big]<\infty.
        \end{align*}
        Particularly, if $\lambda$ is locally Lipschitz with at most linear growth, or if for every $T > 0$
        \begin{align}
        \sup_{t\in[0,T]}\mathbb E[\lambda(P_t)]<\infty. 
        \end{align}
    
    \item[(iv)]
        \textit{Dominance of discount factor.} Under any admissible $(C,\theta)$, moments of $K_t$ and $P_t$ grow at most exponentially with rate strictly smaller than $\rho$.
        
\end{assumptions}

\begin{lemma}
Under the previous assumptions, the transversality condition \eqref{eq:TV} holds.
\end{lemma}

\begin{proof}
Due to the integrability conditions, the result is obtained directly from the growth conditions on $e^{-\rho T}v(K,P)$, the dominance of the discount factor and the dominated convergence theorem.
\end{proof}

\subsection{Viscosity solutions}

\n 

In this section, we develop the corresponding viscosity solution framework. The setting involves a jump-diffusion system with endogenous, pollution-dependent jump intensity and stochastic pollution evolution, leading to an integro-differential HJB equation that may not admit classical smooth solutions. 

In what follows, assume the state process $(K_t,P_t)_{t \geq 0}$ evolves\footnote{Compare with the dynamics from the model presented in Section \ref{section:BM+NHPP-driven-model}.} according to 
\begin{align}
    &dK_t = b^{\mathrm{cap}}(K_t,P_t,C_t,\theta_t)\,dt
	- (1-\omega(K_{t-},P_t))K_{t-}\,dq_t, \label{eq:Visc_K}\\
	&dP_t =  b^{\mathrm{pol}}(K_t,P_t,\theta_t)\,dt + \sigma P_t\,dW_t, \label{eq:Visc_P}
\end{align}
where the drift coefficients are given by \eqref{eq:b_K}--\eqref{eq:b_P}.

Observe that as in the NHPP-driven jump-diffusion model, the infinitesimal generator of the process $(K_t,P_t)_{t \geq 0}$ with respect to the control $(C_t,\theta_t)_{t \geq 0}$ is given by \eqref{eq:Gen-verify-NH}. By applying the DPP over \textit{a short time horizon}, we obtain that $v$ formally satisfies the HJB \eqref{eq:HJB_VerTh} locally. However, since $v$ need not be smooth, we require a more general and robust concept of solution that ensures existence, stability, and uniqueness within the class of continuous functions.

Historically, viscosity solutions provide the appropriate weak notion of solution for such HJB integro-differential equations (see the celebrated reference by Crandall \& Ishii \cite{crandall_users_1992}; alternatively, see \cite{pham_continuous-time_2009} for a more recent reference in the continuous setting, and \cite{ishikawa_optimal_2004} or \cite{barles_second-order_2008} for the discontinuous one). 

Roughly speaking, the idea behind the viscosity approach is to characterize the value function through smooth test functions $\varphi$ that, when evaluated with the integro-differential operator, locally ``touch $v$ from above or below''; that is, they correspond to pointwise (local) maxima or minima of the difference $v - \varphi$. In more rigorous terms, we have the following definition:

\begin{definition}[Viscosity sub- and super-solutions]
	\label{def:ViscosityPham}
	Let $u$ be a continuous function on $\mathcal{S}$.
	\begin{itemize}
		\item[(a)] $u$ is a \emph{viscosity sub-solution} of \eqref{eq:HJB_VerTh} if, for all $\varphi\in C^2(\mathcal{S})$ and all points $(K,P) \in \mathcal{S}$ where $u-\varphi$ attains a local maximum, it holds that
		\begin{equation}\label{eq:visc_sub_Pham}
			\rho u(K,P)
			\le
			\sup_{(C,\theta)\in\mathfrak{a}}
			\big\{
                U(C,P)
    			+ 
                \mathcal L^{C,\theta}\varphi (K,P)
            \big\}.
		\end{equation}
		\item[(b)] $u$ is a \emph{viscosity super-solution} of \eqref{eq:HJB_VerTh} if, for all $\varphi\in C^2(\mathcal{S})$ and all points $(K,P) \in \mathcal{S}$ where $u-\varphi$ attains a local minimum, it holds that
		\begin{equation}\label{eq:visc_super_Pham}
			\rho u(K,P)
			\ge
			\sup_{(C,\theta)\in\mathfrak{a}}
			\big\{
                U(C,P)
    			+ 
                \mathcal L^{C,\theta}\varphi(K,P)
            \big\}.
		\end{equation}
		\item[(c)] $u$ is a \emph{viscosity solution} if it is both a subsolution and a supersolution.
	\end{itemize}
\end{definition}

In order to prove the well-posedness of viscosity solutions in the continuous case, one often relies on the idea of the \textit{sub-} and \textit{super-jets} of a solution $u$, which are suitable generalizations of the notion of the gradient. In the discontinuous case it is still possible -- although not trivial -- to provide a generalization in this direction, despite the fact that the non-locality of the operator $\mathcal{L}^{C,\theta}$ prevents us from taking a \textit{straightforward path}. 

Indeed, in \cite{barles_second-order_2008}, authors surpass this difficulty by providing an equivalent definition of \textit{jets} for which a comparison principle of non-localized test functions hold, along with continuous dependence estimates; in \cite{biswas_viscosity_2010}, a slightly different procedure is taken via a non-local maximum principle of semi-continuous functions; in \cite{imbert_non-local_2005}, a comparison principle is proven for first-order non-local Hamilton-Jacobi equations (i.e. no diffusion term) by means of a regularization of L\'evy kernels. Fortunately, due to the characteristics and structure of our current setting, we can fit many of our models within existent results in the literature. 

The next result establishes that the value function of the control problem is a viscosity solution of the associated HJB equation in the sense of Crandall, Ishii and Lions \cite{crandall_users_1992}, adapted to the integro-differential setting following \cite{barles_second-order_2008}.

\begin{theorem}[Value function as viscosity solution]
	\label{thm:value_viscosity}
	Assume that:
	\begin{enumerate}
		\item[(a)] $U(C,P)$ is continuous and concave in $C$, with polynomial growth;
		\item[(b)] $Y,Z$ are locally Lipschitz and have linear growth;
		\item[(c)] $\omega(K,P)\in(0,1)$ and $\lambda(P)$ are continuous with $\lambda$
		locally Lipschitz;
		\item[(d)] for every admissible control $(C_t,\theta_t)_{t \geq 0}\in\mathcal A(K,P)$,
		the state equations \eqref{eq:Visc_K}--\eqref{eq:Visc_P} admit a unique strong solution.
        
	\end{enumerate}
Then the value function $v$ is continuous and is a viscosity solution of the HJB equation \eqref{eq:HJB_VerTh}.
\end{theorem}

\begin{proof}

In order to prove the result, we show that the function $v$ is continuous on $\cS$ and satisfies both subsolution and supersolution properties.

Fix $T>0$ and define the truncated value function
\[ v_T(K,P) :=
\sup_{(C,\theta)\in\mathcal A(K,P)} \mathbb E \left[\int_0^T e^{-\rho t}U(C_t,P_t)\,dt\right].
\]
Under Assumptions~(b)--(d), for every admissible control $(C,\theta)$ the state process $(K_t,P_t)$ depends continuously in probability on the initial condition $(K,P)$ on compact time intervals. By Assumption~(a) and admissibility, the discounted payoff is uniformly integrable on $[0,T]$. Hence the map
\[
(K,P)\longmapsto
\mathbb E \left[\int_0^T e^{-\rho t}U(C_t,P_t)\,dt\right]
\]
is continuous for each admissible $(C,\theta)$.

Standard $\varepsilon$-optimality arguments imply that $v_T\in C(\mathcal S)$.
Moreover, by discounting and the polynomial growth of $U$, we have
\[
\sup_{(C,\theta)}\mathbb E \left[\int_T^\infty e^{-\rho t}|U(C_t,P_t)|\,dt\right]
\longrightarrow 0 \qquad\text{as }T\to\infty,
\]
uniformly on compact subsets of $\mathcal S$. Hence $v_T\to v$ locally uniformly, and
$v\in C(\mathcal S)$.

Having shown the continuity of $v$, we proceed to verifying its subsolution property.
Let $\varphi\in C^2(\mathcal S)$ and suppose that $v-\varphi$ attains a local maximum
at $(K_0,P_0)\in\mathcal S$, with $v(K_0,P_0)=\varphi(K_0,P_0)$.
Then there exists $r>0$ such that
\[
v(K,P)\le \varphi(K,P)
\quad\text{for all }(K,P)\in B_r(K_0,P_0).
\]
Define the exit time
\[
\tau_r:=\inf\{t\ge0:(K_t,P_t)\notin B_r(K_0,P_0)\},
\qquad
\tau:=\tau_r\wedge h,
\]
where $h>0$.

Fix an arbitrary admissible control $(C,\theta)\in\mathcal A(K_0,P_0)$.
By the dynamic programming principle,
\[
v(K_0,P_0)
\ge
\mathbb E \left[ \int_0^\tau e^{-\rho s}U(C_s,P_s)\,ds + e^{-\rho\tau}v(K_\tau,P_\tau) \right].
\]
Since $(K_\tau,P_\tau)\in\overline{B_r(K_0,P_0)}$ a.s.\ and $v\le\varphi$ on $B_r$,
we obtain
\begin{equation}
\varphi(K_0,P_0) \ge \mathbb E \left[ \int_0^\tau e^{-\rho s}U(C_s,P_s)\,ds + e^{-\rho\tau}\varphi(K_\tau,P_\tau)
\right].
\label{eq:sub1}
\end{equation}

Applying It\^o's formula for jump-diffusions to $e^{-\rho t}\varphi(K_t,P_t)$
on $[0,\tau]$ and taking expectations yields
\begin{align}
\mathbb E \left[e^{-\rho\tau}\varphi(K_\tau,P_\tau)\right]-\varphi(K_0,P_0) =
\mathbb E \left[ \int_0^\tau e^{-\rho s} \big( \mathcal L^{C_s,\theta_s}\varphi(K_s,P_s) - \rho\varphi(K_s,P_s)
\big)\,ds
\right],
\label{eq:sub2}
\end{align}
since the Brownian and compensated jump terms are martingales.

Combining \eqref{eq:sub1} and \eqref{eq:sub2}, we obtain
\begin{equation}
0 \ge \mathbb E \left[ \int_0^\tau e^{-\rho s} \Big( U(C_s,P_s) + \mathcal L^{C_s,\theta_s}\varphi(K_s,P_s) - \rho\varphi(K_s,P_s) \Big)\,ds \right].
\label{eq:sub3}
\end{equation}

Restricting without loss of generality to controls constant on $[0,h]$ and dividing \eqref{eq:sub3} by $h$, dominated convergence and continuity of the coefficients yield,
as $h\downarrow0$,
\[ 0 \ge U(C,P_0) + \mathcal L^{C,\theta}\varphi(K_0,P_0) - \rho\varphi(K_0,P_0).
\]
Since $(C,\theta)\in\mathfrak a$ was arbitrary, we conclude
\[
\rho\varphi(K_0,P_0)
\le
\sup_{(C,\theta)\in\mathfrak a} \Big\{ U(C,P_0) + \mathcal L^{C,\theta}\varphi(K_0,P_0) \Big\},
\]
which is precisely the viscosity subsolution inequality.

Similarly, we show that $v$ satisfies the supersolution property. Let $\varphi\in C^2(\mathcal S)$ and suppose that $v-\varphi$ attains a local minimum
at $(K_0,P_0)$, with equality there.
Fix $\varepsilon>0$. By the dynamic programming principle, there exists an admissible
control $(C^\varepsilon,\theta^\varepsilon)$ such that
\[
v(K_0,P_0) \le \mathbb E \left[
\int_0^\tau e^{-\rho s}U(C^\varepsilon_s,P_s)\,ds + e^{-\rho\tau}v(K_\tau,P_\tau)
\right]+ \varepsilon\,\mathbb E[\tau].
\]
Since $v\ge\varphi$ on $B_r(K_0,P_0)$, we obtain
\begin{equation}
\varphi(K_0,P_0) \le \mathbb E \left[
\int_0^\tau e^{-\rho s}U(C^\varepsilon_s,P_s)\,ds + e^{ \rho\tau}\varphi(K_\tau,P_\tau)
\right]+
\varepsilon\,\mathbb E[\tau].
\label{eq:super1}
\end{equation}

Repeating the It\^o argument above yields
\[ 0 \le \mathbb E \left[ \int_0^\tau e^{-\rho s}
\Big( U(C^\varepsilon_s,P_s) +
\mathcal L^{C^\varepsilon_s,\theta^\varepsilon_s}\varphi(K_s,P_s) - \rho\varphi(K_s,P_s) \Big)\,ds
\right]
+
\varepsilon\,\mathbb E[\tau].
\]
Dividing by $h$, letting $h\downarrow0$, and then $\varepsilon\downarrow0$, we obtain
\[ \rho\varphi(K_0,P_0) \geq \sup_{(C,\theta)\in\mathfrak a} \Big\{
U(C,P_0) + \mathcal L^{C,\theta}\varphi(K_0,P_0) \Big\},
\]
which is the viscosity supersolution inequality.

Thus, $v$ is both a viscosity subsolution and supersolution of the HJB
equation. Together with continuity, this proves that $v$ is a viscosity solution in the sense of Definition~\ref{def:ViscosityPham}.
\end{proof}

\subsubsection{General jump-diffusion model}

\n 

Unlike our previous models, the presence of a random measure with a stochastic intensity kernel $\lambda(P_{t-},\zeta)$ requires some modifications in the verification theorems for viscosity solutions of the value function $v$, due to the stronger interaction between the controlled states and the sizes of the jumps. To overcome this issue, we now make use of the recent developments in infinite-horizon recursive control provided in \cite{luo_infinite_2025}. We focus on a modification of the model \eqref{Eq:K-BM}--\eqref{Eq:K_0-P_0-BM} containing both a Poissonian and a Brownian component and impose an additional set of suitable assumptions:

\begin{assumptions}
\label{Assumptions}
\hfill
\begin{enumerate}
\renewcommand{\theenumi}{\arabic{enumi}}
\renewcommand{\labelenumi}{(\theenumi)}

    \item 
        \label{Assumption-Compactness}
        The set of admissible actions $\mathfrak{a} \subseteq \mathbb{R}_{\geq 0} \times [0,1]$ is compact.

    \item 
        \label{Assumption-Jumps}
        Let $\tilde{\omega} : \mathcal{S} \times \mathbb{R}_{ \geq 0 } \to \mathbb{R}_{\geq 0}$ be a measurable function such that 
        \begin{align*}
            0 \leq \tilde{\omega}( K,P ,\zeta) \leq K, \quad
            &\forall (K, P) \in \mathcal{S}, \quad \forall \zeta \in \mathbb{R}_{\geq 0}.
        \end{align*}
        Furthermore, there exists a measurable mapping $\ell_\gamma : [0, \infty) \to [0,1]$ satisfying
        \begin{align*}
            &| \ell_\gamma |_{ P, 2 } + | \ell_\gamma |_{ P, m } < \infty\
        \end{align*}
        for some $m \geq 2$ with
        \begin{align*}
            &| \ell_\gamma |_{ P, m } := \Big( \int_{[0,\infty)^2} |\ell_\gamma(\zeta,r)|^m\, \nu(d \zeta) \otimes d r \Big)^\frac{1}{m},
            &
            &\forall P \in \mathbb{R}_{> 0},
        \end{align*}
        such that for all $(K,P), (K',P') \in \mathcal{S}$, $\zeta,r \in [0,\infty)$, $(C,\theta) \in \mathfrak{a}$,
        \begin{align*}
            \big| 
                \tilde{\omega}(K,P,\zeta)\, \mathbf{1}_{[0,\lambda(P,\zeta)]}(r) - \tilde{\omega}(K',P',\zeta)\, \mathbf{1}_{[0,\lambda(P',\zeta)]}(r) 
            \big| \hspace{2em}&
            \\
            \leq 
            \ell_\gamma(\zeta, r) \big( | K - K' | + | P - P' | \big).&
        \end{align*}

    \item 
        \label{Assumption-Drift}
        Let $b^{\mathrm{cap}}$ and $b^{\mathrm{pol}}$ be given as in \eqref{eq:b_K} and \eqref{eq:b_P} for some arbitrary functions $Y$ and $Z$, representing the production and efficiency of abatement respectively, such that for some constants $\ell_b \geq 0$, $\alpha_b > 0$ and for all $(K,P), (K',P') \in \cS$, $(C,\theta) \in \mathfrak{a}$,
        \begin{align*}
            |  \tilde b^{\mathrm{cap}}(K,P,C,\theta) - \tilde b^{\mathrm{cap}}(K',P',C,\theta) | + |b^{\mathrm{pol}}(K,P,C,\theta) - b^{\mathrm{pol}}(K',P',C,\theta) | 
            \\
            \leq 
            \ell_b \big( | K - K' | + | P - P' | \big),
            \\
            \big( \tilde b^{\mathrm{cap}}(K,P,C,\theta) - \tilde b^{\mathrm{cap}}(K',P',C,\theta) \big)(K - K') \hspace{10em}&
            \\
            + \big( b^{\mathrm{pol}}(K,P,C,\theta) - b^{\mathrm{pol}}(K',P',C,\theta) \big)(P - P') 
            \\
            \leq -\alpha_b \big( | K - K' |^2 + | P - P'|^2 \big),
        \end{align*}
        where
        \begin{align*}
            \tilde b^{\mathrm{cap}}(K,P,C,\theta) := b^{\mathrm{cap}}(K,P,C,\theta) - \int_{(0,\infty)} \tilde{\omega}(K,P,\zeta)\, \lambda(P,\zeta)\, \nu(d \zeta).
        \end{align*}
    
    \item
        \label{Assumption-Utility}
        The utility function $U(C,P)$ is Lipschitz on $P$.

    \item 
        \label{Assumption-constants}
        Let $m \geq 2$ be as in the previous point. The following inequality holds
        \begin{align*}
            2 \alpha_b - (m-1) \sigma - \frac{2c_m}{m} |\ell_\gamma |_{P,2}^2 - c_m | \ell_\gamma |_{P,m}^m > 0,
        \end{align*}
        with
        \begin{align*}
            c_m := 
            \begin{cases}
                \frac{m(m-1)}{2} - 1,   & \text{if } 2 < m < 3,
                \\
                2m - 4,                 & \text{if } m = 2 \text{ or } m \geq 3.
            \end{cases}
        \end{align*}  
\end{enumerate}
\end{assumptions}

Some comments regarding the previous assumptions: \eqref{Assumption-Compactness} is a well-known  sufficient condition for the minimization of the Hamiltonian at any given instant $t$; Assumptions \eqref{Assumption-Jumps}, \eqref{Assumption-Drift} and \eqref{Assumption-Utility} correspond to standard regularity conditions from the theory of forward-backward SDEs, with \eqref{Assumption-Jumps} being the Lipschitz regularity of the jump component in terms of the extended Poisson space, i.e. the underlying random measure $N$ from Assumptions \ref{Assumption-Poisson-Imbedding}, while also maintaining the economic interpretation that disasters destroy part of the capital; lastly, \eqref{Assumption-constants} is a technical assumption linking the convexity of the value function to the regularity of the jump sizes. This last point is explored with more detail in \cite{luo_infinite_2025}.

\begin{definition}
Let $\mathfrak{a} \subseteq \mathbb{R}_{\geq 0} \times [0,1]$ be compact. For $m \geq 2$ and $t \geq 0$, the set of admissible controls at time $t$ is defined as
\begin{align*}
    \mathcal{A}_{t}^m 
    :=
    \bigg\{
        (C,\theta) : \Omega \times [t, \infty) \to \mathfrak{a}
        ~\bigg|~
        &(C,\theta)\text{ is }\mathbb{F}\text{-predictable, and}
        \\
        &\hspace{3em}\mathbb{E}\Big[ \int_t^\infty \max\{ C_s^2, C_s^m \} ds \Big| \mathcal{F}_t \Big] < \infty
    \bigg\}.
\end{align*}
\end{definition}

With the assumptions set, we are ready to present the main result of this section. In what follows, we use integral notation on the SDEs in order to clarify the forward and backward components.
\begin{theorem}
\label{Thm:Viscosity-with-BSDEs}
Under our standing assumptions, the following hold:
\begin{enumerate}
    \item 
        \label{Thm:Viscosity-with-BSDEs-1}
        For every $(C_t,\theta_t)_{t \geq 0} \in \mathcal{A}_0^m$, there exists a strong unique solution to the system of forward-backward SDEs in infinite-time horizon:
        \begin{align*}
            K_t 
            &= 
            K_0
            +
            \int_0^t \tilde b^{\mathrm{cap}} \big( K_s, P_s, C_s, \theta_s \big)\,dt 
            \\
            &\hspace{4em}- 
            \int_{ (0,t] \times (0,\infty) \times (0,\infty) } \tilde{\omega}\big( K_s, P_s, \zeta \big)\, \mathbf{1}_{[0,\lambda(P_s,\zeta)]}(r) \,\tilde{N}( ds, d\zeta, dr),
            \\
            P_t 
            &= 
            P_0
            +
            \int_0^t b^{\mathrm{pol}} \big( K_s, P_s, C_s, \theta_s \big)\,dt  
            + 
            \int_0^t \sigma P_s\,d W_s,
        \end{align*}
        with $K_0 = K, P_0 = P  > 0$ for all $t \geq 0$, and 
        \begin{align}
            \label{Eq:BSDE}
            \mathcal{V}_t 
            &=
            \mathcal{V}_T
            +
            \int_t^T \big( U(C_s, P_s) - \rho \mathcal{V}_s \big) ds 
            - 
            \int_t^T \mathcal{Z}_s d W_s 
            \\ \nonumber
            &\hspace{4em}-
            \int_{ (t,T] \times (0,\infty) \times (0,\infty) } \mathcal{U}_s(\zeta , r) \tilde{N}( ds, d\zeta, dr),
        \end{align}
        for all $0 \leq t \leq T < \infty$. 

    \item 
        \label{Thm:Viscosity-with-BSDEs-2}
        There exists a unique viscosity solution $v$ to the equation 
        \begin{equation}
            \label{Eq:Brownian-HJB}
            \begin{aligned}
                \rho v(K,P) 
                =
                \sup_{ (C,\theta) \in \mathfrak{a} } \Big\{ 
                    &U(C,P) 
                    + 
                    \frac{\partial v}{\partial K}(K,P) b^{\mathrm{cap}} \big( K,P, C, \theta \big) \\
                    &+
                    \frac{\partial v}{\partial P}(K,P) b^{\mathrm{pol}} \big( K,P, C, \theta \big)
                \Big\} 
                +
                \frac{1}{2} \frac{\partial^2 v}{\partial P^2} (K,P)\,\sigma^2\,P^2 \\
                &+
                \int_{  (0,\infty) } \big( v( \omega(K,P, \zeta )K,P) - v(K,P) \big) \lambda(p,\zeta)\,\nu( d\zeta ).
            \end{aligned}
        \end{equation}
        Moreover, 
        \begin{align}
            \label{Eq:BSDE-initial-value}
            \mathcal{V}_0 
            = 
            \sup_{(C,\theta) \in \mathcal{A}_0^m} J(K,P; C,\theta) 
            = 
            J(K,P; \hat{C},\hat{\theta}) 
            = 
            v(K,P).
        \end{align} 

    \item
        \label{Thm:Viscosity-with-BSDEs-3}
        Conversely, if $v$ is a classical (respectively, viscosity) solution to the HJB \eqref{Eq:Brownian-HJB}, such that for the maximizer of the Hamiltonian at time $t$ there exists an admissible control $(\hat{C}_t, \hat{\theta}_t)_{t \geq s} \in \mathcal{A}_s^m$ for all $s \geq 0$, and the previous forward-backward SDE has a strong solution for the control $(\hat{C}_t, \hat{\theta}_t)_{t \geq 0}$; then $(\hat{C}_t, \hat{\theta}_t)_{t \geq 0} \in \mathcal{A}_0^m$ is an optimal control and $v$ is the corresponding value function. 
        
\end{enumerate}
\end{theorem}

Before presenting the proof of the previous theorem, we briefly illustrate how the BSDE in infinite-time horizon \eqref{Eq:BSDE} encodes the information corresponding to the original objective functional $J$ from \eqref{eq:Optimal_v}. Let $(\mathcal{V},\mathcal{Z},\mathcal{U})$ be a solution to the infinite-horizon BSDE \eqref{Eq:BSDE}. An application of It\^o's formula to $g(t, \mathcal V_t)$ with $g(t,x) = e^{-\rho t} x$ yields
\begin{equation*}
    \begin{aligned}
         e^{-\rho T}\mathcal{V}_T - e^{-\rho t}\mathcal{V}_t 
         =
         &- \int_t^T e^{-\rho s} U( C_s, P_s ) d s 
         + 
         \int_t^T e^{-\rho s} \mathcal{Z}_s d W_s
         \\
         &+
         \int_{ (t,T] \times (0,\infty) \times (0,\infty) } \mathcal{U}_s(\zeta , r) \tilde{N}( ds, d\zeta, dr),
    \end{aligned}
    \qquad \forall 0 \leq t \leq T < \infty.
\end{equation*}
Then, applying conditional expectations with respect to $\mathcal{F}_t$, the following recursive relation for the running utility is obtained:
\begin{align}
    \label{Eq:Recursion}
    &\mathcal{V}_t = \mathbb{E} \Big[ e^{ - \rho ( T - t ) } \mathcal{V}_T + \int_t^T e^{-\rho (s - t)} U( C_s, P_s ) d s \Big| \mathcal{F}_t \Big],
    \qquad \forall 0 \leq t \leq T < \infty.
\end{align}
In other words,
\begin{align}
    \mathcal{V}_t 
    = 
    \liminf_{ T \to \infty } \mathbb{E} \Big[ e^{ - \rho ( T - t ) } \mathcal{V}_T + \int_t^T e^{-\rho (s - t)} U( C_s, P_s ) d s \Big| \mathcal{F}_t \Big],
    \qquad \forall t \geq 0,
\end{align}
from where we can deduce the necessity of the usual transversality condition
\begin{align*}
    \liminf_{ T \to \infty } \mathbb{E} \big[ e^{ - \rho T } \mathcal{V}_T \big] = 0
\end{align*}
in order to obtain
\begin{align*}
    \mathcal{V}_0 = \liminf_{ T \to \infty } \mathbb{E} \Big[ \int_0^T e^{-\rho s} U( C_s, P_s ) d s \Big],
\end{align*}
see Lemma 3.2 in \cite{luo_infinite_2025}. It is worth noting that the type of recursive relation presented in \eqref{Eq:Recursion} -- hence the name \textit{recursive utility} -- has been known in literature; in the realm of randomized settings, we would like to point out to the pioneer work Duffie and Epstein \cite{duffie_stochastic_1992} on stochastic differential utility.

Once the relation between BSDEs and the objective functional has been established, the proof of the main result reduces to verifying the conditions required in \cite{luo_infinite_2025}:

\begin{proof}[Proof of Theorem \ref{Thm:Viscosity-with-BSDEs}]
Let 
\begin{align*}
    f( K, P, C, \theta, \mathcal{V})
    =
    U(C, P) - \rho \mathcal{V},
\end{align*}
and observe that for all $t \geq 0$, $(C,\theta) \in \mathfrak{a}$, $(K,P), (K',P') \in \mathcal{S}$ and $\mathcal{V}, \mathcal{V}' \in \mathbb{R}$,
\begin{align*}
    | f( K, P, C, \theta, \mathcal{V}) - f( K', P', C, \theta, \mathcal{V}') |
    \leq
    \ell_U | P - P' | + \rho | \mathcal{V} - \mathcal{V}' |,
    \\
    \big( f( K, P, C, \theta, \mathcal{V}) - f( K, P, C, \theta, \mathcal{V}') \big) \cdot ( \mathcal{V} - \mathcal{V}' ) 
    \leq 
    - \rho | \mathcal{V} - \mathcal{V}' |^2,
\end{align*}
for some nonnegative Lipschitz constant $\ell_U$. This, coupled with Assumptions \ref{Assumptions}, means that the main conditions $\mathrm{(C1)}_p$, $\mathrm{(C2)}$, $\mathrm{(C3)}$ and $\mathrm{(C4)}'$ of \cite{luo_infinite_2025} hold true, with $p$ and $\rho$ in their paper corresponding to $m$ and $1$ in our current setting, respectively. Consequently, most of the results of Theorem \ref{Thm:Viscosity-with-BSDEs} are obtained from propositions in \cite{luo_infinite_2025}: \ref{Thm:Viscosity-with-BSDEs-1} is due to Lemmas 3.1 and 3.2; the equality of \eqref{Eq:BSDE-initial-value} in \ref{Thm:Viscosity-with-BSDEs-2} corresponds to Theorem 4.7; and the verification results from \ref{Thm:Viscosity-with-BSDEs-3} are obtained through Theorems 5.2 and 5.7 for the classical and viscosity solution cases, respectively.
\end{proof}

Before closing this section, we point out that, although Assumptions \ref{Assumptions} may exclude some of the models considered before, the present framework could be extended to include them. More precisely, $\tilde{\omega}$ and $U$ are now assumed to be \textit{globally} Lipschitz in the state variables, rather than merely \textit{locally} Lipschitz as in previous sections, see Theorem \ref{thm:value_viscosity}. 

This is a technical restriction due to the approach taken by \cite{luo_infinite_2025}, which requires the existence of strong solutions to BSDEs with infinite-time horizon in order to proof the existence of a unique viscosity solution $v$ to the HJB equation \eqref{Eq:Brownian-HJB}. This does not contradict our previous statements, since points \ref{Thm:Viscosity-with-BSDEs-1} and \ref{Thm:Viscosity-with-BSDEs-2} of Theorem \ref{Thm:Viscosity-with-BSDEs} fall within the scope of \ref{thm:value_viscosity}, but it does affect the verification result from the theorem, i.e. point \ref{Thm:Viscosity-with-BSDEs-3}. We \ conjecture that these conditions of global Lipschitz regularity can be weakened, for example, by extending the results in \cite{abdelhadi_locally_2022} that ensure the existence of strong solutions to BSDEs with locally Lipschitz coefficients from a finite to an infinite-time-horizon setting, or the ones found in \cite{papapantoleon_existence_2018} to include locally Lipschitz drivers.

\section{Conclusion}

\n 

This paper develops a unified stochastic control framework for growth-environment models in which the intensity and severity of rare disasters are endogenously linked to the state of pollution.
Building on the setup of \cite{brausmann_escaping_2024}, we formulate an infinite-horizon social planner’s problem with capital and pollution evolving under controlled jump(-diffusion) dynamics, where disasters destroy a state-dependent fraction of capital and may arrive at pollution-dependent rates, possibly with additional marks and sources of randomness. 
Within this framework, we define the value function and characterize optimal trade-offs among consumption, investment, and abatement under environmentally driven catastrophe risk.

From a modeling perspective, we organize several specifications from the economics and climate literature into a single coherent hierarchy. We start with a benchmark model featuring a standard Poisson process with constant intensity, then allow disaster intensity to depend on pollution via a nonhomogeneous Poisson process, capturing feedback from environmental degradation to catastrophe risk. 
Introducing Brownian noise in pollution dynamics leads to jump-diffusion models and integro-differential HJB equations with both local (diffusive) and nonlocal (jump) terms. 
At the most general level, all cases are embedded in a framework based on Poisson random measures with state-dependent compensators and random marks, accommodating random disaster magnitudes and richer environmental dynamics. This nesting clarifies how pollution-driven disaster models extend the constant-intensity benchmark and provides a unified language for comparison.

On the analytical side, we derive the Hamilton-Jacobi-Bellman equations from the dynamic programming principle for pure-jump, jump-diffusion, and Poisson random measure specifications. For sufficiently regular cases, we establish a verification theorem (Theorem~\ref{thm:Pham-verify}) showing that a $C^2$ solution of the HJB equation coincides with the value function and yields optimal feedback controls, and we clarify how the HJB structure and optimality conditions evolve as the jump mechanism becomes richer. 
When such regularity cannot be guaranteed, we show that the value function is a viscosity solution of the HJB equation under natural assumptions (Theorem~\ref{thm:value_viscosity}). 
In the most general setting, with both Brownian and Poisson-driven pollution dynamics, we further characterize the value function via forward-backward stochastic differential equations with jumps. Theorem~\ref{Thm:Viscosity-with-BSDEs} establishes a one-to-one correspondence between solutions of the infinite-horizon FBSDE and classical or viscosity solutions of the HJB equation, thereby linking the stochastic control, its PIDE representation, and the FBSDE approach and providing a flexible analytical and numerical toolkit for studying pollution-driven disaster models.  

Although our focus is on the mathematical structure of pollution-driven disaster models, the flexibility of the random-measure framework allows for a wide range of economically and analytically relevant extensions.
These may include multiple pollutants and capital stocks, threshold-based disaster intensities, and strategic interaction among multiple agents (e.g. regions), leading to game-theoretic control problems. 
On the analytical side, natural extensions involve state-constraint problems (e.g., at zero capital or pollution) and learning about disaster intensity, which give rise to robust control formulations and HJB equations with additional nonlinearities and offer promising avenues for studying climate and environmental policy under deep uncertainty. The introduction of Poisson random measures also allows for further dynamics to be considered such as \textit{Hawkes processes}, a particular type of point process with self-exciting interaction that has been proven relevant in finance \cite{hawkes_hawkes_2022} and nature-related policies \cite{lesage_hawkes_2022}.

\section*{Acknowledgments}
Joshué Helí Ricalde-Guerrero gratefully acknowledges the support of the SNF project MINT 205121-21981.

\bibliographystyle{plain}
\bibliography{references}

\appendix

\section{Technical derivation of HJB equations}
\label{Sec:Technical-appendix}

\medskip

\n 

In this appendix we present the computations needed to obtain the HJB equation for the models presented in the article. We shall follow the same structure as before and divide the results according to the jump-mechanism of the models, i.e. constant-size jumps (Poisson processes over $[0,\infty)$) and random-size jumps (Poisson random measures).

\subsection{Models of constant-size jumps}

\subsubsection{Standard Poisson Process}
\label{section:StandardPP_HJB}

\n

To compute $v(K_h, P_h)$, we apply It\^{o}'s formula for jump processes to the value function $v$, assuming that the state dynamics are given by a jump-diffusion process with jump times driven by Poisson process $(q_t)_{t \geq 0}$ of intensity $\lambda > 0$. Then
\begin{equation} \label{eq:DPP1}
	\begin{aligned} 
		v(K_h, P_{h}) 
		= 
        v(K_0, P_0) 
                &+ 
                \int_{0}^{h} \Big(\frac{\partial v}{\partial K}(K_{s-}, P_{s-}) dK_s 
                + 
                \frac{\partial v}{\partial P}(K_{s-}, P_{s-}) dP_s \Big) 
                \\
                &+ 
                \sum_{s \leq h} \Big(v(K_s, P_s) - v(K_{s-}, P_{s-}) - \frac{\partial v}{\partial K}(K_{s-}, P_{s-}) \Delta K_s\Big),
	\end{aligned}
\end{equation}
where we used that $P_t$ evolves continuously, i.e., $\Delta P_s = 0$, for all $s > 0$, and hence $P_{s-} = P_s$. The jumps in $K_t$ correspond to the discrete losses in capital due to natural disasters.

Substituting \eqref{eq:DPP1} into \eqref{eq:DPP_ineq}, yields
\begin{equation}
	\begin{aligned}
		v(K,P) &\geq \bbE \Bigg[ \int_{0}^{h} e^{-\rho s} U(C_s, P_s)\,ds \\
		&\hspace{2em}+ e^{-\rho h} \bigg( v(K,P) + \int_{0}^{h} \Big(\frac{\partial v}{\partial K}(K_{s-}, P_{s}) dK_s +\frac{\partial v}{\partial P}(K_{s-}, P_{s}) \, dP_s\Big)  \\
		&\hspace{2em}+ \sum_{0 < s \leq h} \Big(v(K_s, P_s) - v(K_{s-}, P_{s-}) - \frac{\partial v}{\partial K}(K_{s-}, P_{s-}) \Delta K_s \Big) \bigg) \Bigg].
	\end{aligned}
\end{equation}

Now, using the state dynamics of capital \eqref{Opt_K} and pollution \eqref{Opt_P}, we obtain

\begin{equation} \label{eq:HJB2}
	\begin{aligned}
		v(K,P) 
        - 
        e^{-\rho h} v(K,P) 
        &\geq 
        \bbE \Bigg[
            \int_{0}^{h} e^{-\rho s} U(C_s, P_s) \,ds 
            \\
            &\hspace{2em}+ 
            \int_{0}^{h} e^{-\rho h} \big[\frac{\partial v}{\partial K}(K_{s-}, P_{s}) \big((1 - \theta_s) Y_s(K_s) - C_s\big)  
            \\
            &\hspace{4em}+ 
            \frac{\partial v}{\partial P}(K_{s-}, P_{s}) \big( (\phi - z  \theta_s) Y_s(K_s) - \alpha P_s\big) \big] \, ds  
            \\
            &\hspace{2em}- 
            \int_{0}^{h} e^{-\rho h} \frac{\partial v}{\partial K}(K_{s-}, P_{s}) (1-\omega(K_{s-}, P_s))K_s \, dq_s  
            \\
            &\hspace{2em}+ 
            \sum_{0 < s \leq h} e^{-\rho h}  \big(v(K_s, P_s) - v(K_{s-}, P_{s}) \big) 
            \\
            &\hspace{2em}- 
            \sum_{0 <s \leq h} e^{-\rho h}  \frac{\partial v}{\partial K}(K_{s-}, P_{s}) \Delta K_s
        \Bigg].
	\end{aligned}
\end{equation}

Jumps in the capital process $(K_t)_{t\ge0}$ occur only at the arrival times of a natural disaster, i.e., at the jump times of the counting process $(q_t)$. At such times, capital is instantaneously reduced to
\begin{equation} \label{eq:capital-jumps}
K_t=\omega(K_{t-},P_{t-})K_{t-},
\end{equation}
where $\omega(K_{t-},P_{t-})\in(0,1)$ denotes the surviving fraction of capital. The corresponding jump size is
\begin{equation}\label{eq:K_jump}
\Delta K_t = -\big(1-\omega(K_{t-},P_{t-})\big)K_{t-}\,\Delta q_t,
\end{equation}
so that capital losses are proportional to the pre-disaster stock.

We now isolate the contribution of capital jumps to the dynamic programming principle. In particular, we focus on the contribution of the jump component to the marginal variation of the value function, that is, the term representing the instantaneous change in the value function induced by the occurrence of a jump. Using \eqref{eq:K_jump}, the jump component in It\^o’s formula yields
\begin{equation}
	\begin{aligned} 
		\sum_{s \leq h} \frac{\partial v}{\partial K}(K_{s-}, P_{s}) \Delta K_s 
		&= - \int_{0}^{h} \frac{\partial v}{\partial K}(K_{s-}, P_{s}) (1 - \omega (K_{s-}, P_s)) K_{s-} \,dq_s,
	\end{aligned}
\end{equation} 
where the integral is understood in the sense of a stochastic integral with respect to the Poisson process.

Substituting this expression into the integral form of the DPP \eqref{eq:HJB2} gives
\begin{equation} \label{eq:HJB3}
	\begin{aligned}
		v(K,P) - e^{-\rho h} v(K,P) & \geq \bbE \bigg[\int_{0}^{h} e^{-\rho s} U(C_s, P_s) \, ds \\
		&\hspace{2em}+ e^{-\rho h} \Big( \int_{0}^{h} \Big[\frac{\partial v}{\partial K}(K_{s-}, P_{s}) \big((1 - \theta_s) Y_s(K_s) - C_s\big) \\
		&\hspace{2em}+ \frac{\partial v}{\partial P}(K_{s-}, P_{s}) \big( (\phi - z  \theta_s) Y_s(K_s) - \alpha P_s\big) \Big] \,ds  \\
		&\hspace{2em}+ \sum_{0 < s \leq h} \big(v(K_s, P_s) - v(K_{s-}, P_{s}) \big) \Big) \bigg].
	\end{aligned}
\end{equation} 

Since capital jumps only at disaster times and satisfies \eqref{eq:capital-jumps}, the jump sum can be written as
\begin{equation} \label{eq:HJB4}
	\begin{aligned}
		\sum_{0 < s \leq h} \big(v(K_s, P_s) - v(K_{s-}, P_{s})\big)
		&= \int_{0}^{h} \big(v(\omega(K_{s-}, P_s) K_{s-}, P_s) - v(K_{s-}, P_{s}) \big) \,dq_s.
	\end{aligned}
\end{equation}  
Substituting \eqref{eq:HJB4} into \eqref{eq:HJB2} yields
\begin{equation} \label{eq:HJB5}
	\begin{aligned}
		v(K,P) - e^{-\rho h} v(K,P) & \geq \bbE \bigg[\int_{0}^{h} e^{-\rho s} U(C_s, P_s) \,ds  \\
		&\hspace{2em}+ e^{-\rho h} \Big( \int_{0}^{h} \Big[\frac{\partial v}{\partial K}(K_{s-}, P_{s}) \big((1 - \theta_s) Y_s(K_s) - C_s \big) \\
		&\hspace{2em}+ \frac{\partial v}{\partial P}(K_{s-}, P_{s}) \big( (\phi - z  \theta_s) Y_s(K_s) - \alpha P_s\big) \Big] \,ds \\
		&\hspace{2em}+ \int_{0}^{h} \big(v(\omega(K_{s-}, P_{s})K_{s-}, P_s) - v(K_{s-}, P_{s}) \big) \,dq_s \Big) \bigg].
	\end{aligned}
\end{equation}  
The expression \eqref{eq:HJB5} represents the key inequality that leads to the HJB equation when dividing by $h$ and letting $h \to 0$:

\begin{equation} \label{eq:HJB6}
	\begin{aligned}
		&\lim\limits_{h \to 0} \frac{v(K,P) - e^{-\rho h} v(K,P)}{h}  \geq \lim\limits_{h \to 0} \bbE \bigg[\int_{0}^{h} e^{-\rho s} U(C_s, P_s) \,ds \bigg] \\
		&\hspace{8em}+ \lim\limits_{h \to 0} \bbE \bigg[ e^{-\rho h} \frac{1}{h}\int_{0}^{h} \frac{\partial v}{\partial K}(K_{s-}, P_{s}) \big((1 - \theta_s) Y(K_s) - C_s\big) \,ds \bigg]\\
		&\hspace{8em}+ \lim\limits_{h \to 0} \bbE \bigg[ e^{-\rho h} \frac{1}{h} \int_{0}^{h} \frac{\partial v}{\partial P}(K_{s-}, P_{s}) \big( (\phi - z  \theta_s) Y(K_s) - \alpha P_s\big) \,ds \bigg] \\
		&\hspace{8em}+ \lim\limits_{h \to 0} \bbE \bigg[ e^{-\rho h} \frac{1}{h} \int_{0}^{h} \big(v(\omega(K_{s-}, P_{s}) K_{s-}, P_s) - v(K_{s-}, P_{s}) \big) \,dq_s \Big) \bigg].
	\end{aligned}
\end{equation} 

To derive the associated HJB equation, we study the infinitesimal generator of the controlled stochastic process. Starting from the dynamic programming inequality for the value function $v(K,P)$ in \eqref{eq:HJB6}, we investigate the limiting behavior of the expression as the time increment $h \to 0$. The analysis proceeds term by term, applying L'H\^opital's rule, the mean value theorem, continuity of the coefficients, and properties of the Poisson process.

Thus, under the standing smoothness assumptions and using standard limit arguments, the individual terms converge as follows: for the left hand side of the inequality, the mean value theorem yields
\begin{align}
    \lim_{h\to0}\frac{v(K,P)-e^{-\rho h}v(K,P)}{h}
    &= 
    \rho v(K,P). \label{eq:HJB_disc} 
\end{align}
For the first term at the right hand side, using Tonelli's theorem and Lebesgue's differentiation theorem, we obtain
\begin{align}
    \lim_{h\to0}\bbE \left[\frac{1}{h}\int_0^h e^{-\rho s}U(C_s,P_s)\,ds\right]
    &= 
    U(C_0,P). \label{eq:HJB_util} 
\end{align}

Next, we consider the contribution of the drift component of capital accumulation. Since the integrand is predictable and continuous in the state variables, and the control processes are admissible, we apply dominated convergence together with the continuity of the coefficients and the value function derivatives. It follows that the time-average of the integrand converges to its value at $t=0$, yielding
\begin{equation}
    \label{eq:HJB_K} 
    \begin{aligned}    
    \lim_{h\to0}\bbE \left[e^{-\rho h}\frac{1}{h}\int_0^h
    \frac{\partial v}{\partial K}(K_{s-},P_s)\big((1-\theta_s)Y(K_s)-C_s\big)\,ds\right]\hspace{2em}
    \\
    = 
    \frac{\partial v}{\partial K}(K,P)\big((1-\theta_0)Y(K)-C_0\big).
    \end{aligned}
\end{equation}
An analogous argument applies to the pollution stock:
\begin{equation}
    \label{eq:HJB_P}
    \begin{aligned}
    \lim_{h\to0}\bbE \left[e^{-\rho h}\frac{1}{h}\int_0^h
    \frac{\partial v}{\partial P}(K_{s-},P_s)\big((\phi-z\theta_s)Y_s(K_s)-\alpha P_s\big)\,ds\right]\hspace{2em}
    \\
    = 
    \frac{\partial v}{\partial P}(K,P)\big((\phi-z\theta_0)Y(K)-\alpha P\big).
    \end{aligned}
\end{equation}
 
The jump term is evaluated using the compensator of the Poisson process via the Doob-Meyer decomposition:

\begin{equation}
dq_s = \tilde{q}_s + \lambda \, ds,
\end{equation}
where $\tilde q_t := q_t-\lambda t$ is the compensated Poisson process. Since $\tilde q$ is an
$\mathbb{F}$-martingale, its integral against any predictable integrand has zero expectation, i.e.,
\begin{equation}
	\bbE \bigg[ \int_{0}^{h} \big(v(\omega(K_{s-}, P_s) K_{s-}, P_s) - v(K_{s-}, P_{s}) \big) \, d\tilde{q}_s \bigg] = 0.
\end{equation}
It follows that
\begin{equation} 
 \begin{aligned}
 	&\lim\limits_{h \to 0} \bbE \left[ e^{-\rho h} \frac{1}{h} \int_{0}^{h} \big(v(\omega(K_{s-}, P_s) K_{s-}, P_s) - v(K_{s-}, P_{s}) \big) \, dq_s \Big) \right]  \\
 	&= \lim\limits_{h \to 0}  \frac{1}{h} \bbE  \left[ \int_{0}^{h} \big(v(\omega(K_{s-}, P_s) K_{s-}, P_s) - v(K_{s-}, P_{s}) \big) \, d\tilde{q}_s \right]  \\
 	& \hspace{3em} + \lim\limits_{h \to 0}  \lambda \bbE  \left[ \frac{1}{h} \int_{0}^{h} \big(v(\omega(K_{s-}, P_s) K_{s-}, P_s) - v(K_{s-}, P_{s}) \,ds \right] \\
    &= \lambda \big(v(\omega(K,P) K,P) - v(K,P) \big).
 \end{aligned}
\end{equation}

Using the obtained results, we arrive at the dynamic programming inequality
\begin{equation}
 \begin{aligned}
 	\rho v(K,P) \geq U(C, P) &+ \frac{\partial v}{\partial P}(K,P) \big( (\phi - z \theta) Y(K) - \alpha P \big) \\
 	&+  \frac{\partial v}{\partial K}(K,P) \big((1 - \theta) Y(K) - C\big) \\
 	&+ \lambda \big(v(\omega(K,P) K,P) - v(K,P) \big).
 \end{aligned}
\end{equation}
Taking the supremum over admissible actions $(C, \theta) \in \mathfrak{a}$, gives the
Hamilton-Jacobi-Bellman equation (see e.g. \cite{garroni_second_2002}, \cite{oksendal_applied_2019} or \cite{pham_continuous-time_2009}).
\begin{equation} \label{eq:HJB_final}
\begin{aligned}
 	\rho v(K,P) = \sup_{(C, \theta) \in \mathfrak{a}} \Big\{U(C, p) &+ \frac{\partial v}{\partial P}(K,P) \big( (\phi - z \theta) Y(K) - \alpha P \big) \\
 	&+ \frac{\partial v}{\partial K}(K,P) \big((1 - \theta) Y(K) - C\big) \\
 	&+ \lambda \big(v(\omega(K,P) K,P) - v(K,P) \big)\Big\}.
\end{aligned}
\end{equation}

The equation \eqref{eq:HJB_final} characterizes the value function $v(K,P)$  as the unique (viscosity) solution (under suitable regularity assumptions) to the associated stochastic optimal control problem involving continuous dynamics and Poisson-driven jump risk.

\subsubsection{Nonhomogeneous Poisson}

The derivation of the dynamic programming principle and the associated HJB equation proceeds in the same way as in the homogeneous case in Section~\ref{section:StandardPP_HJB}, with one key modification: the constant intensity $\lambda$ is replaced by the state-dependent intensity $\lambda(P_{t-})$, with the corresponding compensator being
\begin{equation*} 
    \Lambda_t = \int_0^t \lambda(P_{s-})\,ds,
    \qquad
    \forall t \geq 0.
\end{equation*}

For $h>0$ sufficiently small, applying It\^o’s formula to $v$ and taking expectations yields

\begin{equation} \label{eq:HJB1_NonHP}
	\begin{aligned}
		v(K,P) - e^{-\rho h} v(K,P) & \geq \bbE \bigg[\int_{0}^{h} e^{-\rho s} U(C_s, P_s) \,ds  \\
		&\hspace{2em}+ e^{-\rho h} \Big( \int_{0}^{h} \Big[\frac{\partial v}{\partial K}(K_{s-}, P_{s}) \big((1 - \theta_s) Y_s(K_s) - C_s \big) \\
		&\hspace{2em}+ \frac{\partial v}{\partial P}(K_{s-}, P_{s}) \big( (\phi - z  \theta_s) Y_s(K_s) - \alpha P_s\big) \Big] \,ds \\
		&\hspace{2em}+ \int_{0}^{h} \big(v(\omega(K_{s-},P_s) K_{s-}, P_s) - v(K_{s-}, P_{s}) \big) \,d q_s \Big) \bigg].
	\end{aligned}
\end{equation}  

Using the Doob-Meyer decomposition $d q_s=d \tilde{q}_s+\lambda(P_{s-})\,ds$, the jump term splits into a martingale part and a compensator part:
\begin{equation} \label{eq:jump_split}
	\begin{aligned}
	&\int_0^h \big(v(\omega(K_{s-},P_s) K_{s-}, P_s) - v(K_{s-}, P_s)\big) \, dq_s
	\\
    &= 
    \int_0^h \big(v(\omega(K_{s-},P_s) K_{s-}, P_s) - v(K_{s-}, P_s)\big)\,d\tilde{q}_s 
    \\
	&\hspace{2em}+
    \int_0^h \big(v(\omega(K_{s-},P_s) K_{s-}, P_s) - v(K_{s-}, P_s)\big)\lambda(P_{s-})\,ds.
	\end{aligned}
\end{equation}
Since the martingale integral has zero expectation,
only the compensator contributes:
\begin{equation}
    \begin{aligned}
    \bbE \left[\int_0^h \big(v(\omega(K_{s-},P_s)K_{s-},P_s)-v(K_{s-},P_s)\big)\,d\tilde{q}_s\right]
    =
    0.
    \end{aligned}
\end{equation}

Consequently,
\begin{equation}\label{eq:jump_NHPP_limit}
    \begin{aligned}
    \lim_{h\to0}&\bbE \left[ e^{-\rho h}\frac{1}{h}\int_0^h \big( 
        v(\omega(K_{s-},P_s)K_{s-},P_s)-v(K_{s-},P_s)
    \big)\,d q_s \right] 
    \\
    &= \lambda(P)\big(v(\omega(K,P)K,P)-v(K,P)\big).
    \end{aligned}
\end{equation}
Substituting this result back into \eqref{eq:HJB1_NonHP}, dividing by $h$, and taking the limit as $h \to 0$, we obtain the HJB inequality
\begin{equation}
	\begin{aligned}
		\rho v(K,P) 
        \geq 
        U(C, p) 
            &+ 
            \frac{\partial v}{\partial P}(K,P) \big( (\phi - z \theta) AK - \alpha P \big) 
            \\
            &+
            \frac{\partial v}{\partial K}(K,P) \big((1 - \theta) AK - C\big) 
            \\
            &+ 
            \lambda(P) \big(v(\omega(K,P) K,P) - v(K,P) \big).
	\end{aligned}
\end{equation}
Finally, optimizing over admissible controls $(C, \theta) \in \mathfrak{a}$, gives the Hamilton-Jacobi-Bellman equation for the nonhomogeneous Poisson case:
\begin{equation} \label{eq:HJB_final_NonHP}
	\begin{aligned}
		\rho v(K,P) 
        = 
        \sup_{(C, \theta) \in \mathfrak{a}} \Big\{ 
            U(C, p) 
            &+ 
            \frac{\partial v}{\partial P}(K,P) \big( (\phi - z \theta) AK - \alpha P \big) 
            \\
            &+ 
            \frac{\partial v}{\partial K}(K,P) \big((1 - \theta) AK - C\big) 
            \\
            &+ 
            \lambda(P)\big(v(\omega(K,P) K,P) - v(K,P) \big)\Big\},
	\end{aligned}
\end{equation}
where $\lambda(P) = \lambda_0 + \lambda_1 P$.
The equation \eqref{eq:HJB_final_NonHP} characterizes the value function $v$  as the unique solution (under suitable regularity assumptions outlined before) to the associated stochastic optimal control problem involving continuous dynamics and Poisson-driven jump risk.

\subsubsection{Brownian-driven pollution with nonhomogeneous Poisson jumps}

\n 

Applying It\^o’s formula for jump-diffusions to the value function $v$ over
$[0,h]$ yields
\begin{equation} \label{eq:value_fn_Ito_BM_1}
\begin{aligned}
	v(K_h,P_h)-v(K,P)
	&= \int_0^h \frac{\partial v}{\partial K}(K_{s-},P_{s-})\big((1-\theta_s)Y(K_s)-C_s\big)\,ds \\
	&\quad + \int_0^h \frac{\partial v}{\partial P}(K_{s-},P_{s-})\big(\phi Y(K_s)-Z(\theta_s Y(K_s))-\alpha P_s\big)\,ds \\
	&\quad + \int_0^h \tfrac12 \sigma^2 P_s^2 \frac{\partial^2 v}{\partial P^2}(K_{s-},P_{s-})\,ds + \int_0^h \frac{\partial v}{\partial P}(K_{s-},P_{s-}) \sigma P_s\,dW_s \\
	&\quad + \int_0^h \big(v(\omega(P_s,K_{s-})K_{s-},P_s)-v(K_{s-},P_s)\big)\,d\hat q_s.
\end{aligned}
\end{equation}
Substituting \eqref{eq:Opt_K_NHPP} and \eqref{eq:P_W_NonH_Pois} into \eqref{eq:value_fn_Ito_BM_1}, yields
\begin{equation}\label{eq:value_fn_Ito_BM_2}
	\begin{aligned}
		v(K_h,P_h)-v(K,P)
		&=\int_0^h \frac{\partial v}{\partial K}\big((1-\theta_s)Y(K_s)-C_s\big)\,ds \\
		&\quad+\int_0^h \frac{\partial v}{\partial P}\big(\phi Y(K_s)-Z(\theta_s Y(K_s))-\alpha P_s\big)\,ds \\
		&\quad+\int_0^h \tfrac12 \sigma^2 P_s^2 \frac{\partial^2 v}{\partial P^2}(K_{s-},P_{s})\,ds
		+\int_0^h \frac{\partial v}{\partial P}(K_{s-},P_{s}) \sigma P_s\,dW_s 
        \\
		&\quad+\int_0^h \big(v(\omega(P_s,K_{s-})K_{s-},P_s)-v(K_{s-}\,P_s)\big)\, d\hat q_s.
	\end{aligned}
\end{equation}
Inserting \eqref{eq:value_fn_Ito_BM_2} into \eqref{eq:DPP_ineq} and considering a small time horizon $h > 0$, we obtain:
\begin{equation} \label{eq:value_fn_Ito_BM_E}
	\begin{aligned}
		v(K,P) - e^{-\rho h} v(K,P) 
        & \geq 
        \bbE \bigg[\int_{0}^{h} e^{-\rho s} U(C_s, P_s) \,ds  
        \\
		&\hspace{2em}+ 
        \int_{0}^{h} e^{-\rho h} \Big(\frac{\partial v}{\partial K}(K_{s-}, P_{s}) \big((1 - \theta_s) Y(K_s) - C_s \big) 
        \\
		&\hspace{4em} + 
        \frac{\partial v}{\partial P}(K_{s-}, P_{s}) \big( (\phi - z \theta_s) Y(K_s) - \alpha P_s\big) 
        \\
		&\hspace{4em} + 
        \frac{1}{2} \sigma^2 P^2_s \frac{\partial^2 v}{\partial K^2}(K_{s-}, P_{s}) \Big) \,ds 
        \\
		&\hspace{2em} +  
        e^{-\rho h} \int_{0}^{h} \frac{\partial v}{\partial P}(K_{s-},P_{s})z P_s\, dW_s 
        \\
        &\hspace{2em}+ 
        e^{-\rho h}  \int_{0}^{h} \big(v(\omega(K_{s-},P_s) K_{s-}, P_s) - v(K_{s-}, P_{s}) \big) \,d q_s \bigg].
	\end{aligned}
\end{equation}  
The last term in \eqref{eq:value_fn_Ito_BM_E} represents the contribution of the jump process. 

Using the Doob-Meyer decomposition
$d q_s=d \tilde{q}_s+\lambda(P_{s-})\,ds$, only the compensator contributes:
\begin{equation} \label{eq:BM_martingale_E_1}
\begin{aligned}
&\bbE \left[\int_0^h
\big(v(\omega(K_{s-},P_s)K_{s-},P_s)-v(K_{s-},P_s)\big)\,dq_s
\right] \\
&=
\bbE \left[\int_0^h \big(v(\omega(K_{s-},P_s)K_{s-},P_s)-v(K_{s-},P_s)\big)\lambda(P_{s-})\,ds\right].
\end{aligned}
\end{equation}

The term 
\begin{equation} \label{eq:BM_martingale_2}
\int_0^h e^{-\rho h} \frac{\partial v}{\partial P} \sigma P_s\,dW_s
\end{equation}
is a true martingale on $[0,h]$, and thus
\begin{equation} \label{eq:BM_martingale_E_2}
	\bbE\left[\int_0^h e^{-\rho h} \frac{\partial v}{\partial P}(K_{s-},P_s) \sigma P_s\,dW_s\right] = 0.
\end{equation}

Using \eqref{eq:BM_martingale_E_1} and \eqref{eq:BM_martingale_E_2}, the \eqref{eq:value_fn_Ito_BM_E} becomes
\begin{equation} \label{eq:value_fn_Ito_BM_3}
	\begin{aligned}
		v(K,P) - e^{-\rho h} v(K,P) & \geq \bbE \bigg[\int_{0}^{h} e^{-\rho s} U(C_s, P_s) \,ds  
        \\
		&\hspace{0.69em}+ 
        \int_{0}^{h} e^{-\rho h} \Big(\frac{\partial v}{\partial K}(K_{s-}, P_{s}) \big((1 - \theta_s) Y(K_s) - C_s \big) 
        \\
		&\hspace{3.5em} + 
        \frac{\partial v}{\partial P}(K_{s-}, P_{s}) \big( (\phi - z  \theta_s) Y(K_s) - \alpha P_s\big) 
        \\
		&\hspace{3.5em} +
        \frac{1}{2} \sigma^2 P^2_s \frac{\partial^2 v}{\partial K^2}(K_{s-}, P_{s})\Big) \,ds 
        \\ 
		&\hspace{0.69em}+ 
        \int_{0}^{h} e^{-\rho h}   \big(v(\omega(K_{s-},P_s) K_{s-}, P_s) - v(K_{s-}, P_{s}) \big) \lambda(P_{s-})\,d s \bigg].
	\end{aligned}
\end{equation}

Dividing \eqref{eq:value_fn_Ito_BM_3} by $h$, letting $h \to 0$, and using dominated convergence yields HJB inequality 
\begin{equation} \label{eq:HJB_ineq_NonH_Pois_BM}
	\begin{aligned}
		\rho v(K,P) \geq U(C, P) &+ \frac{\partial v}{\partial P}(K,P) \big( (\phi - z \theta) Y(K) - \alpha P \big) + \frac{1}{2} \sigma^2 P^2 \frac{\partial^2 v}{\partial K^2}(K,P)
        \\
		&+ 
        \frac{\partial v}{\partial K}(K,P) \big((1 - \theta) Y(K) - C\big) \\
        &+ \lambda(P) \big(v(\omega(K,P) K,P) - v(K,P) \big),
	\end{aligned}
\end{equation}
where $\lambda(P)=\lambda_0+\lambda_1 P$.

Optimizing over admissible controls $(C, \theta) \in \mathfrak{a}$, we arrive at the Hamilton-Jacobi-Bellman equation for the Brownian-driven pollution model with nonhomogeneous Poisson jumps
\begin{equation} \label{eq:HJB_NonH_Pois_BM_final}
	\begin{aligned}
		\rho v(K,P) = \sup_{(C, \theta) \in \mathfrak{a}} \Big\{U(C, P) 
        &+ 
        \frac{\partial v}{\partial P}(K,P) \big( (\phi - z \theta) Y(k)  - \alpha P \big) 
        + 
        \frac{1}{2} \sigma^2 P^2 \frac{\partial^2 v}{\partial K^2}(K,P) \\
        &+  
        \frac{\partial v}{\partial K}(K,P) \big((1 - \theta) Y(K)  - C\big)\Big\} \\
		&+ 
        \lambda(P) \big(v(\omega(K,P) K,P) - v(K,P) \big).
	\end{aligned}
\end{equation}

\subsection{Models of random-size jumps}

\n 

We now derive the HJB equation for the model introduced in Section \ref{Eq:PRM-Model-Introduction}. The argument proceeds in two steps: we first obtain the corresponding equation for the intermediate model from Section \ref{Section:An_intermediate_model} (i.e. with the Gaussian component removed), and then we extend the result to the randomized pollution model.

\subsubsection{Equation for the intermediate model}

\n 

The overall strategy remains the same, up to some minor modifications necessitated by $q$ being a Poisson random measure. First, note that instead of \eqref{Eq:Delta_q}, size of the jumps vary according to $\Delta$:
\begin{align*}
    \Delta q_t (d\zeta) 
    = 
    \lim_{h \downarrow 0} q( (t - h, t], d\zeta )
    =
    q \big( \{ t \}, d\zeta \big) = \sum_{n \geq 1} \mathbf{1}_{ \{ \tau_n = t \} } \delta_{\{\Delta_n\}}(d \zeta),
\end{align*}
where the law of $\Delta_n$ is determined as in \eqref{Eq:Conditional-mark-distribution}. Without loss of generality, assume $K_t = K_0$ for all $t < 0$; then, the jump in capital \eqref{eq:K_jump} is replaced by
\begin{equation}
    \label{Eq:Delta_K}
    \begin{aligned}
        \Delta K_t 
        = 
        K_t - K_{t-} 
        &=
        - \lim_{h \downarrow 0} \int_{ (t - h,t] \times [0,\infty) } \big( 1 - \omega( K_{s-}, P_s, \zeta ) \big) K_{s-}\,q( ds, d\zeta)
        \\
        &=
        - \int_{ (0,\infty) } \big( 1 - \omega( K_{t-}, P_t, \zeta ) \big) K_{t-}\,q( \{ t \}, d\zeta)
        \\
        &=
        \begin{cases}
            - \big( 1 - \omega( K_{t-}, P_t, \Delta_n ) \big) K_{t-}
                & \text{on the event }\{ \tau_n = t \},
            \\
            0
                & \text{in any other case,}
        \end{cases}
    \end{aligned}
\end{equation}
for every $t \geq 0$.

With these considerations in mind, applying It\^o's rule (see Chapter 14 in \cite{cohen_stochastic_2015}) on $v(K_t,P_t)$ for the prescribed control $(C_t,\theta_t)$ yields the following equality
\begin{align*}
    v(K_t, P_t) - v(K,P) 
    =& \int_0^t \frac{\partial v}{\partial K} ( K_{s-}, P_{s} ) d K_s + \int_0^t \frac{\partial v}{\partial P} ( K_{s-}, P_{s} ) d P_s 
    \\
    &+
    \sum_{ 0 < s \leq t } \Big( 
        v(K_s, P_s) - v(K_{s-}, P_{s})
        -
        \frac{\partial v}{\partial K} (K_{s-}, P_{s}) \Delta K_s 
    \Big).
\end{align*}
From the dynamic programming principle (see e.g. \cite{oksendal_applied_2019} or \cite{luo_infinite_2025}), equation \eqref{eq:DPP_ineq} holds, and as a consequence we get that for any small $h > 0$,
\begin{equation}
    \label{Eq:DPP-inequality}
    \begin{aligned}
        \frac{v(K,P) - e^{-\rho h} v(K,P) }{h}
        &\geq
        \frac{1}{h} \mathbb{E} \bigg[ 
            \int_0^h e^{-\rho s} U(C_s, P_s) \, ds 
            + 
            e^{-\rho h} \int_0^h \frac{\partial v}{\partial K} ( K_{s-}, P_{s} )\,d K_s 
            \\
            &\hspace{4em}+
            e^{-\rho h} \int_0^h \frac{\partial v}{\partial P} ( K_{s-}, P_{s} )\,d P_s 
            \\
            &\hspace{4em}+
            e^{-\rho h} \sum_{ 0 < s \leq h } v(K_s, P_s) - v(K_{s-}, P_{s})
            \\
            &\hspace{4em}-
            e^{-\rho h} \sum_{ 0 < s \leq h } \frac{\partial v}{\partial K} (K_{s-}, P_{s}) \Delta K_s 
        \bigg].
    \end{aligned}
\end{equation}

We now analyze each term individually when $h$ goes to zero. First, using the same arguments as above, 
we obtain that
\begin{align}
    \label{Eq:Estimate-1}
    \lim_{h \to 0} \frac{v(K,P) - e^{-\rho h} v(K,P)}{h} &= \rho v(K,P),
    \\ \label{Eq:Estimate-2}
    \lim_{h \to 0} \mathbb{E} \bigg[ \frac{1}{h}\int_0^h e^{-\rho s} U(C_s, P_s) \, ds \bigg] &= U(C_0,P),
\end{align}
and
\begin{align}
    \label{Eq:Estimate-3}
    &\lim_{h \to 0} e^{-\rho h}~\mathbb{E} \bigg[ 
        \frac{1}{h} \int_{0}^{h} \frac{\partial v}{\partial P}(K_{s-}, P_{s})\,d P_s 
    \bigg] 
    =
    \frac{\partial v}{\partial P}(K,P) b^{\mathrm{pol}} \big( K,P, C, \theta \big),
\end{align}
respectively.

For the integral with respect to the capital, observe that \eqref{Eq:K} can be rewritten in terms of the compensated martingale measure $\tilde{q}$:
\begin{align*}
    dK_t 
    =& 
    \bigg( 
        b^{\mathrm{cap}} \big( K_t, P_t, C_t, \theta_t \big)
        -
        \int_{ (0,\infty) } \big( 1 - \omega( K_{t-}, P_t, z ) \big) K_{t-} \lambda(P_t, \zeta)\,\nu( d\zeta)
    \bigg)\,dt 
    \\
    &- 
    \int_{ (0,\,\cdot\,]\times (0,\infty) } \big( 1 - \omega( K_{t-}, P_t, \zeta ) \big) K_{t-}\,\tilde{q}( dt, d\zeta),
\end{align*}
where
\begin{align*}
    \tilde{q}( dt, d\zeta)
    :=
    q( dt, d\zeta ) - \Lambda( dt, d\zeta ).
\end{align*}
Then,
\begin{align*}
    &\mathbb{E} \bigg[ \int_0^h \frac{\partial v}{\partial K} ( K_{s-}, P_{s} )\,d K_s \bigg]
    \\
    &=
    \mathbb{E} \bigg[ \int_0^h \frac{\partial v}{\partial K} ( K_{s-}, P_{s} ) b^{\mathrm{cap}} \big( K_s, P_s, C_s, \theta_s \big)\,ds \bigg]
    \\
    &\quad -
    \mathbb{E} \bigg[ \int_{ (0,h] \times (0,\infty) } 
        \frac{\partial v}{\partial K} ( K_{s-}, P_{s} ) \big( 1 - \omega( K_{s-}, P_s, \zeta ) \big) K_{s-} 
    \,\Lambda ( ds, d\zeta ) \bigg].
\end{align*}
On the one hand, from \eqref{Eq:K} we have
\begin{equation}
    \label{Eq:Estimate-4}
    \begin{aligned}
        &\lim_{h \to 0} e^{-\rho h} \mathbb{E} \bigg[ 
            \frac{1}{h} \int_{0}^{h} \frac{\partial v}{\partial K}(K_{s-}, P_{s}) b^{\mathrm{cap}} \big( K_s, P_s, C_s, \theta_s \big)
        \,ds \bigg] \\
        &=
        \frac{\partial v}{\partial K}(K,P) b^{\mathrm{cap}} \big( K,P, C_0, \theta_0 \big).
    \end{aligned}
\end{equation}
On the other hand, from \eqref{Eq:lambda}, \eqref{Eq:Delta_K} and the definition of integral with respect to Poisson random measures,
\begin{equation}
    \label{Eq:Estimate-5}
    \begin{aligned}
        &\mathbb{E} \bigg[ \sum_{ 0 < s \leq h } \frac{\partial v}{\partial K} (K_{s-}, P_{s}) \Delta K_s \bigg]
        \\
        &=
        -\mathbb{E} \bigg[ \int_{ (0,h] \times (0,\infty) } 
            \frac{\partial v}{\partial K} ( K_{s-}, P_{s} ) \big( 1 - \omega( K_{s-}, P_s, \zeta ) \big) K_{s-} 
        \,\Lambda( ds, d\zeta ) \bigg],
    \end{aligned}
\end{equation}
where the last equality is due to Campbell's theorem and the definition of compensator, see \cite{bremaud_point_2020}.

Lastly, observe that
\begin{align*}
    &\mathbb{E} \bigg[ \sum_{ 0 < s \leq h } v(K_s, P_s) - v(K_{s-}, P_{s}) \bigg]
    =
    \mathbb{E} \bigg[ \sum_{ 0 < s \leq h } v(\Delta K_{s} + K_{s-}, P_s) - v(K_{s-}, P_{s}) \bigg]
    \\
    &=
    \mathbb{E} \bigg[ \sum_{0 < s \leq h } \sum_{n \geq 1} 
        \big( v( \omega( K_{s-}, P_s, \zeta_n )K_{s-}, P_s) - v(K_{s-}, P_{s}) \big) 
    \mathbf{1}_{ \{ \tau_n = s \} } \bigg]
    \\
    &=
    \mathbb{E} \bigg[ \int_{ (0, h] \times (0,\infty) }
        \big( v( \omega( K_{s-}, P_s, \zeta )K_{s-}, P_s) - v(K_{s-}, P_{s}) \big)
    \,q( ds, d\zeta ) \bigg]
    \\
    &=
    \mathbb{E} \bigg[ \int_{ (0, h] \times (0,\infty) }
        \big( v( \omega(K_{s-},P_s)( K_{s-}, P_s, \zeta )K_{s-}, P_s) - v(K_{s-}, P_{s}) \big) 
    \,\Lambda( ds, d\zeta ) \bigg],
\end{align*}
where we have again used Campbell's theorem. Then, by Lebesgue differentiation theorem \cite{folland_real_1999},
\begin{equation}
    \label{Eq:Estimate-6}
    \begin{aligned}
        &\lim_{h \downarrow 0} \frac{1}{h}\mathbb{E} \bigg[ \int_{ (0, h] \times (0,\infty) }
            \big( v( \omega( K_{s-}, P_s, \zeta )K_{s-}, P_s) - v(K_{s-}, P_{s}) \big) 
        \,\Lambda( ds, d\zeta ) \bigg]
        \\
        &\hspace{4em}= \int_{  (0,\infty) } \big( v( \omega( K,P, \zeta )K,P) - v( K,P ) \big) \lambda(P,\zeta)\,\nu( d\zeta ).
        \\
        \implies &\lim_{h \downarrow 0} \frac{e^{-\rho h}}{h}\mathbb{E} \bigg[ \int_{ (0, h] \times (0,\infty) }
            \big( v( \omega( K_{s-}, P_s, \zeta )K_{s-}, P_s) - v(K_{s-}, P_{s}) \big) 
        \,\Lambda( ds, d\zeta ) \bigg]
        \\
        &\hspace{4em}= 
        \lim_{h \downarrow 0} e^{-\rho h} 
            \cdot \int_{  (0,\infty) } \big( v( \omega( K,P, \zeta )K,P) - v( K,P ) \big) \lambda(P,\zeta)\,\nu( d\zeta )
        \\
        &\hspace{4em}= 
        1 \cdot \int_{  (0,\infty) } \big( v( \omega( K,P, \zeta )K,P) - v( K,P ) \big) \lambda(P,\zeta)\,\nu( d\zeta ).
    \end{aligned}
\end{equation}
Taking limits at both sides of \eqref{Eq:DPP-inequality} and plugging in the estimates from \eqref{Eq:Estimate-1} to \eqref{Eq:Estimate-6},
we obtain that for any admissible control such that $(C_0,\theta_0) = (C,\theta)$,
\begin{align*}
    \rho v(K,P) 
    \geq&~
    U(C,p) 
    + 
    \frac{\partial v}{\partial K}(K,P) b^{\mathrm{cap}} \big( K,P, C, \theta \big)
    +
    \frac{\partial v}{\partial P}(K,P) b^{\mathrm{pol}} \big( K,P, C, \theta \big)
    \\
    &+
    \int_{  (0,\infty) } \big( v( \omega( K,P, \zeta )K,P) - v( K,P ) \big) \lambda(P,\zeta)\,\nu( d\zeta ).
\end{align*}
Thus, under suitable regularity conditions, the value function $v$ solves 
the Hamilton-Jacobi-Bellman partial integro-differential equation
\begin{equation}
    \label{Eq:Poisson-Measure-HJB}
    \begin{aligned}
        \rho v(K,P) 
        =
        \sup_{ (C,\theta) \in \mathfrak{a} } \Big\{ 
            U(C,P) 
            &+ 
            \frac{\partial v}{\partial K}(K,P) b^{\mathrm{cap}} \big( K,P, C, \theta \big) \\
            &+
            \frac{\partial v}{\partial P}(K,P) b^{\mathrm{pol}} \big( K,P, C, \theta \big) \\
            &+
            \int_{  (0,\infty) } \big( v( \omega( K,P, \zeta )K,P) - v( K,P ) \big) \lambda(P,\zeta)\,\nu( d\zeta ) 
        \Big\}.
    \end{aligned}
\end{equation}

\subsubsection{Extension to randomized pollution}
\label{Section:Model-with-BM+PRM}

\n 

We now present another extension of the model to include a randomized pollution:
\begin{align}
    \tag{\ref{Eq:K-BM}}
	dK_t 
        &= 
        b^{\mathrm{cap}}\big( K_t, P_t, C_t, \theta_t \big)\,dt 
        - 
        \int_{ [0,\infty) } \big( 1 - \omega( K_{t-}, P_t, \zeta ) \big) K_{t-}\,q( dt, d\zeta),
    \\
    \tag{\ref{Eq:P-BM}}
    dP_t 
        &= 
        b^{\mathrm{pol}} \big( K_t, P_t, C_t, \theta_t \big)\,dt + \sigma P_t\,d W_t,
    \\
    \tag{\ref{Eq:K_0-P_0-BM}}
    K_0 & > 0,\quad P_0  > 0,
\end{align}
for some given constant $\sigma > 0$.

We now apply It\^o's rule \cite{cohen_stochastic_2015} on $e^{-\rho t}v(K_t,P_t)$ for the prescribed control $(C_t,\theta_t)$ on the new dynamics:
\begin{align*}
    &e^{-\rho t}v(K_t, P_t) - v(K,P) 
    =
    -\int_0^t \rho e^{-\rho s} v( K_{s-}, P_{s-} ) d s
    \\
    &\hspace{1em}+
    \int_0^t e^{-\rho s} \bigg\{ 
        \frac{\partial v}{\partial K} ( K_{s-}, P_{s-} ) d K_s + \frac{\partial v}{\partial P} ( K_{s-}, P_{s-} ) d P_s 
    \bigg\}
    \\
    &\hspace{1em}+
    \frac{1}{2} \int_0^t e^{-\rho s} \bigg( 
        \frac{\partial^2 v}{\partial K^2 } ( K_{s-}, P_{s-} ) d \langle K^c, K^c \rangle_s 
        +
        \frac{\partial^2 v}{\partial K \partial P} ( K_{s-}, P_{s-} ) d \langle K^c, P^c \rangle_s
        \\
        &\hspace{6em}+
        \frac{\partial^2 v}{\partial P \partial K} ( K_{s-}, P_{s-} ) d \langle P^c, K^c \rangle_s
        +
        \frac{\partial^2 v}{\partial P^2} ( K_{s-}, P_{s-} ) d \langle P^c, P^c \rangle_s
    \bigg)
    \\
    &\hspace{1em}+
    \sum_{ 0 < s \leq t } e^{-\rho s} \bigg( 
        v(K_s, P_s) - v(K_{s-}, P_{s-})
        \\
        &\hspace{6em}-
        \frac{\partial v}{\partial K} (K_{s-}, P_{s-}) \Delta K_s 
        -
        \frac{\partial v}{\partial P} (K_{s-}, P_{s-}) \Delta P_s
    \bigg),
\end{align*}
where $K^c$ and $P^c$ denote the continuous components of $K$ and $P$, respectively. Then, using the same arguments as before and including the Brownian component of $P$, we have that
\begin{align*}
    &e^{-\rho t}v(K_t, P_t) - v(K,P) 
    =
    -\int_0^t \rho e^{-\rho s} v( K_{s-}, P_{s} ) d s
    \\
    &+
    \int_0^t e^{-\rho s} \Big( 
        \frac{\partial v}{\partial K} ( K_{s-}, P_{s} ) d K_s + \frac{\partial v}{\partial P} ( K_{s-}, P_{s} ) d P_s 
    \Big)
    \\
    &\hspace{5em}+
    \frac{1}{2} \int_0^t e^{-\rho s} \frac{\partial^2 v}{\partial P^2} ( K_{s-}, P_{s} ) d \langle P^c, P^c \rangle_s
    \\
    &\hspace{5em}+
    \sum_{ 0 < s \leq t } e^{-\rho s} \Big( 
        v(K_s, P_s) - v(K_{s-}, P_{s})
        -
        \frac{\partial v}{\partial K} (K_{s-}, P_{s-}) \Delta K_s 
    \Big).
\end{align*}
As in the previous cases, from the dynamic programming principle equation \eqref{eq:DPP_ineq} holds, and as a result
\begin{align*}
    &\frac{1}{h} \mathbb{E} \Big[ \rho \int_0^h e^{-\rho s} v( K_{s-}, P_{s} ) \,d s - v(K,P) \Big]
    \\
    &\geq 
    \frac{1}{h}\mathbb{E} \left[ \int_0^h e^{-\rho s} U(C_s, P_s) \, ds  \right]
    +
    \frac{1}{h} \mathbb{E} \Big[ \int_0^h e^{-\rho s} \frac{\partial v}{\partial K} ( K_{s-}, P_{s} ) d K_s \Big]
    \\
    &\hspace{1em}+
    \frac{1}{h} \mathbb{E} \Big[ \int_0^h e^{-\rho s} \frac{\partial v}{\partial P} ( K_{s-}, P_{s} ) d P_s \Big]
    +
    \frac{1}{2h} \mathbb{E} \Big[ \int_0^h e^{-\rho s} \frac{\partial^2 v}{\partial P^2} ( K_{s-}, P_{s} ) d \langle P^c, P^c \rangle_s \Big]
    \\
    &\hspace{1em}+
    \frac{1}{h} \mathbb{E} \Big[ \sum_{ 0 < s \leq t } e^{-\rho s} \Big( 
        v(K_s, P_s) - v(K_{s-}, P_{s})
        -
        \frac{\partial v}{\partial K} (K_{s-}, P_{s-}) \Delta K_s 
    \Big) \Big].
\end{align*}

Observe that the only difference from the previous case is in the inclusion of a Brownian motion in the dynamics of $P$ and in the integral with respect to $\langle P^c \rangle$; however, from the properties of $W$ we have that
\begin{align*}
    &\lim_{h \downarrow 0}\frac{1}{h} \mathbb{E} \Big[ \int_0^h e^{-\rho s} \frac{\partial v}{\partial P} ( K_{s-}, P_{s} ) d P_s \Big]
    \\
    &\hspace{3em}=
    \lim_{h \downarrow 0}\frac{1}{h} \mathbb{E} \Big[ 
        \int_0^h e^{-\rho s} \frac{\partial v}{\partial P} ( K_{s-}, P_{s} ) b^{\mathrm{pol}} \big( K_s, P_s, C_s, \theta_s \big)\,ds
    \Big]
    \\
    &\hspace{3em}=
    \frac{\partial v}{\partial P}(K,P) b^{\mathrm{pol}} \big( K,P, C, \theta \big),
\end{align*}
and
\begin{align*}
    &\lim_{h \downarrow 0} \frac{1}{2h} \mathbb{E} \Big[ 
        \int_0^h e^{-\rho s} \frac{\partial^2 v}{\partial P^2} ( K_{s-}, P_{s} )\,d \langle P^c, P^c \rangle_s 
    \Big]
    =
    \frac{1}{2} \frac{\partial^2 v}{\partial P^2} ( K,P ) \sigma^2 P^2.
\end{align*}
Adding this estimates to the ones presented in the previous section, i.e. equations \eqref{Eq:Estimate-1} through \eqref{Eq:Estimate-6}, yields the following HJB equation:
\begin{equation}
    \label{Eq:Brownian-HJB-Appendix}
    \begin{aligned}
        \rho v(K,P) 
        =
        \sup_{ (C,\theta) \in \mathfrak{a} } \Big\{ 
            U(C,P) 
            &+ 
            \frac{\partial v}{\partial K}(K,P) b^{\mathrm{cap}} \big( K, P, C, \theta \big) \\
            &+
            \frac{\partial v}{\partial P}(K,P) b^{\mathrm{pol}}\big( K, P, C, \theta \big)
        \Big\}
        +
        \frac{1}{2} \frac{\partial^2 v}{\partial P^2} ( K,P ) \sigma^2 P^2\\
        &+
        \int_{  (0,\infty) } \big( v( \omega( K, P, \zeta )K,P) - v( K, P ) \big) \lambda(P,\zeta)\,\nu( d\zeta ).
    \end{aligned}
\end{equation}

\end{document}